\documentclass[11pt,twoside]{amsart}

\usepackage{amssymb}
\usepackage{latexsym}

\newcommand{\Id}{\mathrm{Id}}
\newcommand{\eins}{\mathbf{1}}
\newcommand{\N}{{\mathbb N}}
\newcommand{\dopu}{{:}\allowbreak\ }
\newcommand{\eps}{\varepsilon}
\newcommand{\cal}{\mathcal}
\newcommand{\R}{{\mathbb R}}
\newcommand{\MM}{{\cal M}}
\newcommand{\DD}{{\cal D}}

\newcommand{\mytilde}[1]{\mathbin{\tilde#1}}

\newcommand{\loglike}[1]{\mathop{\mathrm {#1}}\nolimits}

\newcommand{\dist}{\loglike{dist}}
\newcommand{\supp}{\loglike{supp}}
\newcommand{\SDt}{\loglike{{\cal S}\!{\cal D}}}
\newcommand{\narr}{\loglike{{\cal N}\!\!{\cal A}\!{\cal R}}}
\newcommand{\Op}{\loglike{{\cal O}\!{\cal P}}}
\newcommand{\SRN}{\loglike{{\cal S}\!{\cal R}\!{\cal N}}}
\newcommand{\UB}{\loglike{{\cal U}\!{\cal B}}}
\newcommand{\MUB}{\loglike{{\cal M}\!{\cal U}\!{\cal B}}}
\newcommand{\Fin}{{{\cal F}}}
\newcommand{\Lin}{\loglike{lin}}
\newcommand{\conv}{\loglike{conv}}
\newcommand{\clconv}{\loglike{\overline{conv}}}
\newcommand{\cp}{\loglike{cp}}

\newcommand{\rest}[2]{#1\raisebox{-0.3ex}{\mbox{$\mid_{#2}$}}}

\theoremstyle{plain}
\newtheorem{thm}{Theorem}[section]

\newtheorem{theo}[thm]{Theorem}
\newtheorem{prop}[thm]{Proposition}
\newtheorem{cor}[thm]{Corollary}
\newtheorem{lemma}[thm]{Lemma}

\theoremstyle{definition}
\newtheorem{definition}[thm]{Definition}

\theoremstyle{remark}
\newtheorem{example}[thm]{Example}
\newtheorem{rem}[thm]{Remark}

\numberwithin{equation}{section}

\newcounter{abc}   
\newcounter{iiiii} 

\newenvironment{aequivalenz}
{\setcounter{iiiii}{0}
\begin{list}%
{{\rm (\roman{iiiii})}}
{\usecounter{iiiii}
\parsep=0pt plus 1pt
\topsep=1pt plus 2pt minus 1pt
\itemsep=1pt plus 2pt minus 1pt
\leftmargin=3\baselineskip
\labelsep=.6\baselineskip
\labelwidth=2.4\baselineskip
\rightmargin 0pt}%
}%
{\end{list}}

\newenvironment{statements}%
{\setcounter{abc}{0}
\begin{list}%
{{\rm (\alph{abc})}}
{\usecounter{abc}
\parsep=0pt plus 1pt
\topsep=1pt plus 2pt minus 1pt
\itemsep=1pt plus 2pt minus 1pt
\leftmargin=3\baselineskip
\labelsep=.6\baselineskip
\labelwidth=2.4\baselineskip
\rightmargin 0pt}%
}%
{\end{list}}

\begin{document}

\title[Narrow operators and rich subspaces]{Narrow operators and 
rich subspaces of Banach spaces with the Daugavet property}

\author{Vladimir M. Kadets, Roman V. Shvidkoy and Dirk Werner}

\address{Faculty of Mechanics and Mathematics, Kharkov National University,
\qquad {}\linebreak pl.~Svobody~4,  61077~Kharkov, Ukraine}
\email{vishnyakova@ilt.kharkov.ua}

\curraddr{Department of Mathematics, Freie Universit\"at Berlin,
\mbox{Arnimallee~2--6}, D-14\,195~Berlin, Germany}
\email{kadets@math.fu-berlin.de}

\address{Department of Mathematics, University of Missouri,
Columbia MO 65211}
\email{shvidkoy\_r@yahoo.com}

\address{Department of Mathematics, Freie Universit\"at Berlin,
Arnimallee~2--6, \qquad {}\linebreak D-14\,195~Berlin, Germany}
\email{werner@math.fu-berlin.de}

\thanks{The work of the first-named author
was supported by a grant from the {\it Alexander-von-Humboldt
Gesellschaft}.}




\begin{abstract}
Let $X$ be a Banach space. We introduce a formal approach which
seems to be useful in the study of those properties of operators
on $X$ which depend only on the norms of the images of elements.
This approach is  applied to the Daugavet equation for norms
of operators; in particular we develop a general theory of narrow
operators and rich subspaces of spaces $X$ with the Daugavet property
previously studied in the context
of the classical spaces $C(K)$ and $L_{1}(\mu)$.
\end{abstract}

\maketitle

\thispagestyle{empty}

\section{Introduction}

Following \cite{KadSSW} we say that a Banach space $X$ has the
{\em Daugavet property}\/ if for every operator $T\dopu X\to X$ of rank~$1$
the {\em Daugavet equation}
\begin{equation}  \label{eqDE}
\|\Id + T\| = 1 + \|T\|
\end{equation}
is fulfilled. It is known that then every weakly compact operator,
even every strong Radon-Nikod\'ym operator, and every operator
not fixing a copy of $\ell_{1}$ satisfies (\ref{eqDE}) as well
(\cite{KadSSW}, \cite{Shv1}). Incidentally, this shows that our
definition of the Daugavet property is equivalent to the ones which have been
proposed in \cite{AbraAB} and \cite{Kadets}. Classical results due
to Daugavet \cite{Daug}, Lozanovskii \cite{Loz}, and Foia\c{s}, Singer and
Pe\l czy\'nski \cite{FoiSin} state that $C(K)$, $L_{1}(\mu)$ and
$L_\infty(\mu)$  have the Daugavet property provided that
$K$ is perfect and $\mu$ is non-atomic. Recently, corresponding
results in the non-commutative setting were obtained by Oikhberg~\cite{Oik}.

The papers \cite{KadSSW} and \cite{Shv1} study Banach spaces
with the Daugavet property from a structural point of view; for example it
is shown that such a space never embeds into a space having
an unconditional basis, and it contains (many) subspaces isomorphic
to $\ell_{1}$. Also, hereditary properties of the Daugavet property
are established there.

Returning to the classical spaces $C(K)$ and $L_{1}(\mu)$
we mention that a
 different approach to (\ref{eqDE}) on these spaces
 was launched earlier in \cite{KadPop}
and \cite{PliPop}.  These papers study a duality between certain
operators, called {\em narrow operators}, and certain  subspaces,
called {\em rich subspaces}, of such spaces. (For the definitions,
which differ in the two cases, see
Section~\ref{sec2}.) One of the key features of this approach is
that the concept of a narrow operator on $C(K)$ or $L_{1}(\mu)$,
which makes sense for operators from these spaces into an arbitrary
range space, only depends on the values $\|Tx\|$, but not on the
images $Tx$ themselves.

The idea of the present paper is to introduce narrow operators and rich
subspaces in general. In Section~\ref{sec2} we propose a formalism in
order to deal with those properties of an operator which depend only
on the norms of the images of elements.
We define corresponding equivalence
classes and their formal sums and differences, which is reminiscent of
certain procedures in the theory of operator ideals. Then, in
Section~\ref{sec3} we introduce and study narrow operators on Banach
spaces with the Daugavet property. We show, in particular, that strong
Radon-Nikod\'ym operators are narrow and
that narrow operators mapping $X$ to itself satisfy (\ref{eqDE}). In
Section~\ref{sec4} we prove that operators not fixing a copy of
$\ell_{1}$ are narrow, thus extending a result from \cite{Shv1}. To
do so we need an extension of a theorem due to Rosenthal
characterising separable Banach spaces that fail to contain isomorphic
copies of $\ell_{1}$ (Theorem~\ref{ROS1}), which seems to be of
independent interest.

Section~\ref{sec5} deals with rich subspaces. As in the classical
case of $C(K)$ or $L_1(\mu)$, 
a closed subspace $Y\subset X$ is called rich if the quotient
map $q\dopu X\to X/Y$ is narrow. One of the main results here is that
the Daugavet property passes to rich subspaces, which leads to new hereditary
properties.
We also study a related class of
subspaces which we term wealthy. What looks like a quibble of words
is another main result from Section~\ref{sec5}: a subspace is rich if
and only if it is wealthy.
In fact, we also need to deal with a slightly more general class
of operators called strong Daugavet operators.
It turns out that there are strong Daugavet
operators which are not narrow; an example to this effect is
presented in Section~\ref{sec6}.

As for notation, we denote the closed unit ball of a Banach space by
$B(X)$ and its unit sphere by $S(X)$. The slice of $B(X)$ determined
by a functional $x^*\in S(X^*)$ and $\eps>0$ is the set
$$
S(x^*,\eps)= \{x\in B(X)\dopu x^*(x)\ge 1-\eps \}.
$$

We shall repeatedly make use of the following characterisation of the
Daugavet property in terms of slices or weakly open sets from
\cite[Lemma~2.2]{KadSSW} and \cite[Lemma~2.2]{Shv1} respectively.

\begin{lemma} \label{lem1.1}
The following assertions are equivalent:
\begin{aequivalenz}
\item
$X$ has the Daugavet property.
\item
For every $x\in S(X)$, $\eps>0$ and every slice $S$
of $B(X)$ there exists some $y\in S$ such that $\|x+y\|>2-\eps$.
\item
For every $x\in S(X)$, $\eps>0$ and every nonvoid relatively weakly
open subset $U$
of $B(X)$ there exists some $y\in U$ such that $\|x+y\|>2-\eps$.
\end{aequivalenz}
\end{lemma}

Actually, this lemma characterises {\em Daugavet pairs}\/ $(Y,X)$, meaning
a Banach space $X$ and a subspace $Y\subset X$ such that
$$
\|J+T\|=1+\|T\|
$$
for every operator from $Y$ into $X$ of rank~$1$; here $J$ denotes
the canonical embedding map. The only modification to be made in the
formulation of Lemma~\ref{lem1.1} is that $S$ and $U$ refer to slices
and subsets of $B(Y)$.

Finally we mention that all the Banach spaces in this paper are
tacitly assumed to be real.

\bigskip\noindent
\textbf{Acknowledgement.}
We are grateful to the referee whose detailed suggestions have led to a
substantial improvement of the exposition of this paper.


\section{The semigroup $\Op(X)$} \label{sec2}

Throughout the paper the symbol $X$ will  be used for a fixed  Banach space,
the symbols $T$, $ T_i$ etc.\  for bounded  linear operators, acting
from  $X$ to some other  Banach space (not necessarily  the same one).

\begin{definition}
We say that two operators  $ T_1$  and $ T_2$ are {\em equivalent}\/
(in symbols $ T_1 \sim  T_2$)
if $\|T_1x\|=\|T_2x\|$ for every $x\in X$. A class ${\cal M}$ of
operators is said to be {\em admissible}\/
if for every $T \in {\cal M}$ all the members of
the equivalence class of  $T$ also belong to ${\cal M}$.
\end{definition}

In other words, the operators $ T_1$  and $ T_2$ are equivalent if
there is an isometry $U\dopu  T_1(X) \rightarrow  T_2(X)$ such that
$T_2 = U T_1$. For example, the classes of finite-rank operators,
compact  operators,
weakly  compact  operators, operators bounded from below are
admissible; surjections, isomorphisms, projections are examples
of non-admissible operator classes.

\begin{definition}
We say that $ T_1 \le  T_2$ if $\|T_1x\| \le \|T_2x\|$ for
every $x\in X$. A class ${\cal M}$ of
operators forms an {\em order ideal}\/ if for every $T \in {\cal M}$ every
operator  $ T_1 \le  T$  also belongs to ${\cal M}$.
\end{definition}

In other words, $ T_1 \le  T_2$   if
there is a bounded operator  $U\dopu  T_2(X) \rightarrow  T_1(X)$ of norm
${\le1}$ such that
$T_1 = U T_2$. 
Order ideals are clearly admissible.
The classes of finite-rank operators, compact  operators,
weakly  compact  operators are order ideals.

\begin{definition}
A sequence $(T_n)$ of operators is said to be {\em $\sim$convergent}\/
to an operator
$T$ if  $\|T_nx\| \rightarrow  \|Tx\|$ uniformly on $B(X)$. In
terms of $\sim$convergence we define  the notions of
a $\sim$closed set of operators, $\sim$closure, etc.\
in a natural way.
\end{definition}

Of course, the $\sim$limit of a sequence is not unique, but it is unique up
to equivalence of operators.

For example, the  class $\Fin(X)$ of  finite-rank operators
on an infinite-dimensional space $X$ is not $\sim$closed:
its $\sim$closure contains all  compact  operators.
Indeed, let $T\dopu X\to Y$ be compact. Then, for the canonical
isometry $U$ from $Y$ into $C(B(Y^*))$, $T_{1}:=UT$ is compact, too,
and by definition $T_{1}\sim T$. Since $C(B(Y^*))$ has the
approximation property, $T_{1}$ can be approximated by finite-rank
operators in the above sense.

In fact, the $\sim$closure of $\Fin(X)$
coincides with the class ${\cal C}(X)$ of  all  compact
operators since ${\cal C}(X)$ is $\sim$closed. 
To see this suppose that $(T_{n})$ is a $\sim$convergent
sequence of compact operators
on $X$ with limit~$T$. Let $(x_{n})$ be a bounded sequence
in $X$; using a diagonal procedure one can find a subsequence
$(x_{n}')$ such that $(T_{k}x_{n}') ^{\vphantom{\prime}}_{n}$
is convergent for each~$k$.
But $\|T_{k}x\|\to \|Tx\|$ uniformly on bounded sets as $k\to\infty$;
hence $(Tx_{n}')$ is a Cauchy sequence and thus convergent.

\begin{definition}
Let ${\cal N}$ be a  collection of subsets in $X$. We define a
class of operators ${\cal N}^{\sim}$ as follows:
$T \in {\cal N}^{\sim}$ if for every $A \in {\cal N}$, $T$ is
unbounded from below on $A$; i.e.,
\[
\forall \eps>0\ \exists x\in A\dopu \ \|Tx\|\le\eps.
\]
\end{definition}

Evidently, ${\cal N}^\sim$ is a $\sim$closed order ideal, and 
it is {\em homogeneous}\/ in the sense that $\lambda T\in {\cal N}^\sim$
whenever $\lambda\in \R$ and $ T\in {\cal N}^\sim$.
For example, if
${\cal N}=\{S(X)\}$, then ${\cal N}^\sim = \UB(X)$, the class of
operators that are {\em unbounded from below}\/ which is defined by
$$
T\in \UB(X) \iff \inf\{\|Tx\|\dopu \|x\|=1 \} =0.
$$

A significant example for us is the class of all {\em $C$-narrow
operators}\/ on the space $C(K)$. This class was introduced in
\cite{KadPop} as the class of those operators $T\dopu C(K)\to Y$
which are unbounded from below on the unit sphere of each subspace
$J_{F}:= \{f\in C(K)\dopu \rest{f}{F}=0\}$, where $F$ is a proper
closed subset of $K$. To put it another way, if $\cal N$ denotes
the collection of these unit spheres, then the class of $C$-narrow
operators is just ${\cal N}^\sim$.

Another important example is the class of all {\em $L_1$-narrow
operators}\/ on the space $L_1=L_1(\Omega, \Sigma, \mu)$. An
operator $T\dopu L_{1}\to Y$ is called $L_{1}$-narrow if for each
$B\in\Sigma$ and each $\eps>0$ there is a function vanishing off
$B$ and taking only the values $-1$ and $1$ on $B$
such that $\|Tf\|\le\eps$. In other words,   if $\cal N$ now
denotes the collection of these sets of functions, then the class of
$L_{1}$-narrow operators is again ${\cal N}^\sim$. $L_{1}$-narrow
operators were formally introduced in \cite{PliPop}, but the
complement of this class was studied earlier by Ghoussoub and
Rosenthal who called non-$L_{1}$-narrow operators {\em
norm-sign-preserving}. An operator is not $L_{1}$-narrow if and
only if it is not a sign-embedding on any $L_{1}(B)$-subspace
(\cite{GhoRos}, \cite{Ros-signemb}, \cite{Ros-emb}).

We caution the reader that in \cite{KadPop} and \cite{PliPop}
only
the term ``narrow'' is used. In this paper we prefer to speak of
$C$- and $L_{1}$-narrow operators in order not to mix up these
notions with our concept of a narrow operator in
Section~\ref{sec3}; cf., however, Theorem~\ref{cor3.6}.

\bigskip

We now define $\Op(X)$ as the class of all operators on $X$ with the
convention that equivalent operators will be identified. Hence
$\Op(X)$ is actually a collection of equivalence classes, and in fact
it is a set. Namely, for an operator $T$ on $X$ its
equivalent class can be identified with the seminorm $x\mapsto
\|Tx\|$, and the collection of seminorms on $X$ is clearly a set.
Thus, admissible families of operators can be identified with subsets
of $\Op(X)$, and it makes sense to write $T\in\Op(X)$ or ${\cal
M}\subset \Op(X)$.

We now introduce addition and subtraction on $\Op(X)$.
If $T_{1}\dopu X \to Y_{1} $ and
$T_{2}\dopu X \to Y_{2} $ are two operators, define
$$
T_{1}\mytilde+ T_{2}\dopu X \to Y_{1}\oplus_{1}Y_{2}, \quad
x \mapsto (T_{1}x, T_{2}x);
$$
i.e.,
$$
\|(T_{1} \mytilde+ T_{2})x\| = \|T_{1}x\| + \|T_{2}x\|.
$$

\begin{definition}
If ${\cal M}_{1}, {\cal M}_{2} \subset \Op(X)$ are non-empty, then their
{\em $\sim$sum}\/ is defined by ${\cal M}_{1} \mytilde+ {\cal M}_{2}
= \{ T_{1} \mytilde+ T_{2} \dopu T_{1}\in {\cal M}_{1},\ 
T_{2}\in {\cal M}_{2} \}$.
      Their {\em $\sim$difference}\/ is defined by
${\cal M}_{2} \mytilde- {\cal M}_{1}
= \{T\in \Op(X)\dopu T \mytilde+ T_{1} \in {\cal M}_{2}$ whenever
$T_{1} \in {\cal M}_{1} \}$.
\end{definition}

The  operation $ \mytilde +$ is a commutative and associative
operation on  ${\Op}(X)$, and we have $0\in {\cal M}_{2} \mytilde-
{\cal M}_{1}$ if and only if $ {\cal M}_{1} \subset {\cal M}_{2}$.

Let us give some examples.

\begin{example} \label{exam2.7}
Let $K$ be a compact Hausdorff space and
let $\MUB(C(K))$ denote the  class of operators equivalent to some
multiplication operator $M_{h}\dopu f\mapsto hf$ on $C(K)$
which is unbounded
from below; i.e., where $h$ has a zero. Then
 ${\UB}(C(K)) \mytilde - \MUB(C(K))$ consists exactly of the
$C$-narrow operators described above.
\end{example}

\begin{proof}
Let $T\dopu C(K)\to Y$ be $C$-narrow. If $h$ has a zero, we have to
show that, given $\eps>0$,
there is some $f\in S(C(K))$ such that both $\|Tf\|\le\eps$
and $\|hf\|_{\infty}\le\eps$. Now, if $F=\{|h| \ge\eps\}$, which is a
proper subset of $K$, and $f\in
S(J_{F})$ such that $\|Tf\|\le\eps$, then $\|hf\|_{\infty}\le \eps$
as well.

Conversely, if a closed proper subset $F\subset K$ is given, pick
some $h\in S(C(K))$ such that $h=1$ on $F$, $h=0$ off a neighbourhood
$V$ of $F$. If $\|f\|_{\infty}\le1$, $\|Tf\|\le\eps$ and
$\|hf\|_{\infty}\le\eps$, then in
particular $|f|\le\eps$ on $F$. Hence it is possible to replace $f$
by a function $g\in S(J_{F})$ such that $\|Tg\|\le2\eps$, which proves
that $T$ is $C$-narrow.
\end{proof}

For our next example
recall that an operator on $X$ is a left semi-Fredholm operator if
its kernel is finite-dimensional and its range is closed, and it is
strictly singular if it is unbounded from below on (the unit sphere
of) each infinite-dimensional subspace of $X$.

\begin{example}
The class
${\UB}(X)\mytilde -{\Fin}(X)$ consists of all  operators that are not
 left semi-Fredholm operators;
${\UB}(X)\mytilde -({\UB}(X)\mytilde -\Fin(X))$
consists of all strictly singular operators.
\end{example}

\begin{proof}
Let us denote  ${\cal G}(X) = \UB(X) \mytilde- {\cal F}(X)$  
and
${\cal H}(X) = \UB(X) \mytilde- {\cal G}(X)$.

If $T$ is a left semi-Fredholm operator, then, since $\ker T$ is
complemented by a finite-codimensional subspace $Y\subset X$,
$\rest{T}{Y}$  is bounded from below, because $T$ acts as an
isomorphism from $Y$ onto $T(X)$. On the other hand, if $\rest{T}{Y}$
is bounded from below on some finite-codimensional subspace
$Y\subset X$, then $T(Y)$ and, consequently, $T(X)$ must be closed,
and $\ker T$ is finite-dimensional, since otherwise $Y \cap \ker T
\neq \{0\}$. This shows that $T$ is not a left semi-Fredholm operator
if and only if
\begin{statements}
\item
$\rest{T}{Y}$ is not bounded from below on any finite-codimensional
subspace $Y\subset X$.
\end{statements}

Now, if $T$ satisfies (a), $F\in {\cal F}(X)$ and $Y=\ker F$, then
$T\mytilde+ F\in \UB(X)$, i.e., $T\in {\cal G}(X)$. Conversely, if
$T\in {\cal G}(X)$, $Y\subset X$ is finite-codimensional
and $q\dopu X\to X/Y$ is the quotient map,
then, since $T\mytilde+ q\in \UB(X)$, $T$ satisfies~(a).

Thus, we have shown the announced characterisation of ${\cal G}(X)$
and, moreover, we have shown that (a) provides another
characterisation of ${\cal G}(X)$. It follows from (a) that $T\in
{\cal G}(X)$ if and only if
\begin{statements}
\item[(b)]
for every $\eps>0$ there exists an infinite-dimensional subspace
$Z\subset X$ such that $\|\rest{T}{Z}\|\le \eps$;
\end{statements}
see \cite[Prop.~2.c.4]{LiTz1}.

 From (b) it is clear that every strictly singular operator belongs
to ${\cal H}(X)$. On the other hand, if $S$ is not strictly singular
and is bounded from below on some infinite-dimensional subspace $Z$,
then we have for the quotient map $q\dopu X\to X/Z$ that $S\mytilde+
q$ is bounded from below. Since $q$ obviously satisfies (b), this
shows that $S\notin {\cal H}(X)$.
\end{proof}

Let us list some elementary properties of the operation $\mytilde -$
that follow directly from the definition.

\begin{prop}\label{prop2.7}
Suppose that $  {\cal M}_1, {\cal M}_2 \subset \Op(X)$ 
contain  the zero operator.
\begin{statements}
\item
$ {\cal M}_2 \mytilde -  {\cal M}_1 $ is an order ideal or
$\sim$closed whenever ${\cal M}_{2}$ is.
\item
If ${\cal M}_{1}$ and ${\cal M}_{2}$ are order ideals, then 
$ {\cal M}_2 \mytilde -  {\cal M}_1 $ is homogeneous whenever 
$ {\cal M}_2 $ is.
\end{statements}
\end{prop}

Of particular relevance are subsets of $\Op(X)$ that are semigroups
with respect to the operation~$\mytilde+$. 

\begin{prop}\label{prop2.8}
Suppose that $  {\cal M}_1, {\cal M}_2 \subset \Op(X)$ 
contain  the zero operator.
\begin{statements}
\item 
${\cal M}_1$ is a subsemigroup of ${\Op}(X)$ if and only if
${\cal M}_1 \mytilde - {\cal M}_1 \supset {\cal M}_1$, in which case
${\cal M}_1 \mytilde - {\cal M}_1 = {\cal M}_1$.
\item
Let ${\cal M}_1$ be a subsemigroup of ${\Op}(X)$, and let
 ${\cal M}_1 \subset {\cal M}_{2}$. Then
${\cal M}_{2} \mytilde - ({\cal M}_{2} \mytilde -  {\cal M}_1)$ is again
a subsemigroup.
\item 
${\cal M}_{2} \mytilde -  {\cal M}_2$ is always a subsemigroup of
$\Op(X)$.
\end{statements}
\end{prop}

\begin{proof}
(a) is clear from the definition. 

For (b) we note first that 
\begin{equation} \label{eq2.1}
\MM_{2}\mytilde- (\MM_{2} \mytilde- ( \MM_{2}-\MM_{1} )) = \MM_{2}-\MM_{1} .
\end{equation}
Indeed, by definition of $\mytilde-$ we have
\begin{equation} \label{eq2.2}
\MM_{2}\mytilde- (\MM_{2} \mytilde- \MM_{1}) \supset \MM_{1} ,
\end{equation}
whence
$$
\MM_2 \mytilde- ( \MM_{2}\mytilde- (\MM_{2} \mytilde- \MM_{1})) 
\subset \MM_{2} \mytilde- \MM_{1} .
$$
On the other hand, an
application of (\ref{eq2.2}) with $\MM_{1}$ replaced by 
$\MM_{2} \mytilde- \MM_{1}$ gives
``$\supset$'' in (\ref{eq2.1}). Now, by elementary arithmetic
involving $\mytilde+$ and $\mytilde-$ we have, writing $\DD = 
\MM_{2} \mytilde- \MM_{1}$ for short,
\begin{align*}
(\MM_{2} \mytilde- \DD) \mytilde-  (\MM_{2} \mytilde- \DD)
&=
\MM_{2} \mytilde- (\DD \mytilde+ (\MM_{2} \mytilde- \DD))
\displaybreak[0]\\
&=
\MM_{2} \mytilde- ( (\MM_{2} \mytilde- \DD) \mytilde+ \DD)
\displaybreak[0]\\
&=
(\MM_{2}        \mytilde-  (\MM_{2} \mytilde- \DD) ) \mytilde- \DD
\displaybreak[0]\\
&=
\DD \mytilde- \DD \qquad\text{(by (\ref{eq2.1}))}
\displaybreak[0]\\
&=
(\MM_{2} \mytilde- \MM_{1}) \mytilde- \DD
\\
&=
\MM_{2} \mytilde- (\MM_{1} \mytilde+ \DD).
\end{align*}
Because $\MM_{1}$ is a semigroup, one can easily deduce that
$\MM_{1} \mytilde+ \DD \subset \DD$; 
indeed,
\begin{align*}
\MM_{1} \mytilde+ \DD 
&=
( \MM_{2} \mytilde- \MM_{1} ) \mytilde+ \MM_{1}
\displaybreak[0]\\
&=
( \MM_{2} \mytilde- ( \MM_{1} \mytilde+ \MM_{1} )) \mytilde+ \MM_{1}
\displaybreak[0]\\
&=
( ( \MM_{2} \mytilde-  \MM_{1} ) \mytilde- \MM_{1} ) \mytilde+ \MM_{1}
\\
&\subset
\MM_{2} \mytilde- \MM_{1}.
\end{align*}
Therefore
$$
(\MM_{2} \mytilde- \DD) \mytilde-  (\MM_{2} \mytilde- \DD)
\supset \MM_{2} \mytilde- \DD,
$$
completing the proof that 
$ \MM_{2} \mytilde- ( \MM_{2}   \mytilde- \MM_{1} )$ is a semigroup.

Finally, (c) is the special case ${\cal M}_{1} = \{0\}$ of~(b).
\end{proof}

The following definition is important for our abstract semigroup
approach.

\begin{definition}
 Let ${\cal M} \subset {\Op}(X)$, and let ${\cal M}_1 \subset
 {\cal M}$ be  a subsemigroup of ${\Op}(X)$.  ${\cal M}_1$ is
called a {\em maximal subsemigroup}\/ of ${\cal M}$ if every subsemigroup
${\cal M}_2 \subset  {\cal M}$ which includes ${\cal M}_1$
coincides with ${\cal M}_1$.
We call the intersection of all maximal subsemigroups of ${\cal M}$
the {\em central part}\/ of ${\cal M}$ and denote it by $\cp({\cal M})$.
\end{definition}

Here is a characterisation of the central part of $\cal M$.

\begin{thm} \label{cpart1}
Let ${\cal M} \subset {\Op}(X)$ have the following
properties: $0 \in {\cal M}$ and every element of
${\cal M}$ is contained in a
subsemigroup of ${\cal M}$ (this happens for example if ${\cal M}$
is homogeneous). Then $\cp({\cal M})= {\cal M} \mytilde - {\cal M}$.
\end{thm}

\begin{proof}
Let ${\cal M}_1$ be a maximal subsemigroup of ${\cal M}$.
Put ${\cal M}_2={\cal M} \mytilde - {\cal M}$. We have proved above
in Proposition~\ref{prop2.8}(c) that ${\cal M}_2$ is a subsemigroup, so
${\cal M}_2 \mytilde+ {\cal   M}_1$  is a subsemigroup, too. By
definition of  ${\cal M}_2$ we have
${\cal M}_2 \mytilde+ {\cal   M}_1 \subset {\cal M}$. So the
maximality of ${\cal M}_1$ implies that $ {\cal   M}_1 \supset {\cal M}_2$.
This proves the inclusion $\cp({\cal M})\supset {\cal M} \mytilde - {\cal M}$.

Let us now prove the inverse inclusion. Let
$T \in \cp({\cal M})\setminus ({\cal M} \mytilde - {\cal M})$.
Then there is some $T_1 \in {\cal M}$ such that $T_1 \mytilde+ T$
does not belong to ${\cal M}$.
Consider the maximal subsemigroup ${\cal M}_3$
of ${\cal M}$ which contains $T_1$. Then ${\cal M}_3$
cannot contain $T$, so $\cp({\cal M})$  cannot contain $T$ either.
\end{proof}

For every operator $T$ and $\varepsilon > 0$  we define the \textit{tube}
$$
U_{T,\varepsilon}=\{x \in X\dopu  \|Tx\|< \varepsilon\}.
$$
Let ${\cal M}\subset \Op(X)$.
 Put
$$
{\cal M}_{\sim}= \{U_{T,\varepsilon}\cap S(X)\dopu T \in {\cal M},\
\varepsilon > 0 \}
$$
Then  $({\cal M}_{\sim})^{\sim} ={\UB}(X)\mytilde - {\cal M}$.

\begin{prop}\label{prop1}
Let ${\cal M}\subset \Op(X)$ and let  ${\cal N}$ be a collection of
subsets in $X$. Then ${\cal N}^{\sim} \mytilde - {\cal M}=
{\cal N_1}^{\sim}$, where ${\cal N_1}$ consists of all
intersections  of the form $U_{T,\varepsilon} \cap A$,
$T \in {\cal  M}$,  $A \in {\cal N}$, $\varepsilon > 0$.
In particular, if ${\cal N}^{\sim} \mytilde - {\cal M}$ is
non-empty, then all the
intersections   $U_{T,\varepsilon} \cap A$ are non-empty
and ${\cal N}^{\sim} \supset {\cal M}$.
\end{prop}

\begin{proof}
Let $T_1 \in {\cal N}^{\sim} \mytilde - {\cal M}$. Then
for every $T \in {\cal  M}$ we have $T_1 \mytilde + T \in {\cal N}^{\sim}$.
This means that for every $A \in {\cal N}$ and $\varepsilon > 0$
there is an element $x \in A$ such that
$\|(T_1 \mytilde + T)x\| < \varepsilon$. This in turn implies that
$x \in A \cap U_{T,\varepsilon}$ and $\|T_1 x\| < \varepsilon$.
So $T_1 \in {\cal N_1}^{\sim}$.

Now let $T_1 \in {\cal N_1}^{\sim}$. Then
for every $T \in {\cal  M}$, every $A \in {\cal N}$ and $\varepsilon > 0$
there is an element $x \in A \cap U_{T,\varepsilon/2}$ such that
$\|T_1 x\| < \varepsilon/2$. But by the definition of tubes,
$\|T x\| < \varepsilon/2$. So $\|(T_1 \mytilde + T)x\| < \varepsilon$
and $T_1 \in {\cal N}^{\sim} \mytilde - {\cal M}$.
\end{proof}


\section{Narrow operators} \label{sec3}

In this section we define the class of narrow operators on a Banach
space with the Daugavet property. But first we need to introduce a
closely related class of operators.

\begin{definition} \label{defi3.1}
An  operator $T\in\Op(X)$ is said to be a {\em strong Daugavet operator}\/
if for every two elements $x, y \in S(X)$ and  for every
$\varepsilon > 0$ there is an element
$z \in (y + U_{T,\varepsilon}) \cap  S(X)$ such that
$\|z+x\|> 2 - \varepsilon$.
We denote the class of all strong Daugavet operators on $X$
by $\SDt(X)$.
\end{definition}

It follows from Lemma~\ref{lem1.1} that a finite-rank operator on a
space with the Daugavet property
 is a strong Daugavet operator, and conversely, if every rank-$1$
operator is strongly Daugavet, then $X$ has the Daugavet property.

There is an obvious connection between strong Daugavet operators and
the Daugavet equation.

\begin{lemma}
If $T\dopu X\to X$ is a strong Daugavet operator, then $T$ satisfies
the Daugavet equation {\rm(\ref{eqDE})}.
\end{lemma}

\begin{proof}
We assume without loss of generality that $\|T\|=1$.
Given $\eps>0$ pick $y\in S(X)$ such that $\|Ty\|\ge
1-\eps$. If $x= Ty/\|Ty\|$ and $z$ is chosen according to
Definition~\ref{defi3.1}, then
$$
2-\eps< \|z+x\| \le \|z+Ty\| + \eps \le \|z+Tz\| + 2\eps,
$$
hence
$$
\|z+Tz\|\ge 2-3\eps ,
$$
which proves the lemma.
\end{proof}

We now relate the strong Daugavet property to a collection of subsets
of $X$.

\begin{definition}
For every ordered pair of elements $(x, y)$ of $S(X)$
and every $\varepsilon > 0$ let us define a
set $D(x ,y, \varepsilon)$  by
$$
z  \in D(x ,y, \varepsilon) \iff 
\|z+x+y\|> 2 - \varepsilon ~\&~ \|z+y\|< 1 + \varepsilon.
$$
By ${\cal D}(X)$ we denote the collection of all  sets
$D(x ,y, \varepsilon)$, where  $x, y \in S(X)$
and $\varepsilon > 0$.
\end{definition}

\begin{prop}  \label{theo3.3}
 $\SDt(X) =  {\cal D}(X)^{\sim}$.
\end{prop}

\begin{proof}
$T \in {\cal D}(X)^{\sim}$ if and only if for every pair
$x, y \in S(X)$ and $\varepsilon > 0$ there is an element
$z  \in D(x ,y, \varepsilon)$ such that $ \|Tz\|< \varepsilon$.
This in turn is equivalent to the following condition: for every pair
$x, y \in S(X)$ and $\varepsilon > 0$ there is an element $v$
such that $\|v\|< 1 + \varepsilon$, $\|x+v\|>2-\varepsilon$
and $v$ belongs to the
tube $ y + U_{T,\varepsilon}$ (just put $v=z+y$). Evidently,
the last equation coincides with the strong Daugavet property of the operator
$T$.
\end{proof}

Let us consider an example.

\begin{thm} \label{theo3.4}
For a compact Hausdorff space $K$,
 the class $\SDt(C(K))$ of strong Daugavet operators
coincides with the  class of $C$-narrow operators.
\end{thm}

\begin{proof}
The fact that every $C$-narrow operator is  a
strong Daugavet operator has been proved in a  slightly different
form in  \cite[Th.~3.2]{KadPop}. Consider the converse implication. Let
$T \in \SDt(C(K))$. Fix a closed subset $F\subset K$
and  $0<\varepsilon < 1/4$. According to the definition it is sufficient
to prove that  there is a function $f \in S(C(K))$
for which the  restriction to
$F$ is less than $2 \varepsilon$
and $\|Tf\| < 2 \varepsilon$ (cf.\ Example~\ref{exam2.7}).
Let us fix a neighbourhood $U$ of $F$ and an open set $V \subset K$, $V\cap
U=\emptyset$.
Select inductively functions $x_n, y_n  \in S(C(K))$
and $f_{n}, g_n  \in C(K)$
as follows.
All the  $y_n$ are supported on
$ U$, and  the $x_n$ are non-negative functions
supported on $ V$. Given $x_{n}$ and $y_{n}$ pick $f_n \in D(x_n ,y_n,
\varepsilon)$ with
$\|Tf_n\| < \varepsilon$, and let
$g_n = f_1+ f_2 + \cdots + f_n$. Then choose $y_{n+1}\in S(C(K))$
subject to the above support condition such that
$ \sup_{t\in F} |g_n(t)| \, y_{n+1} $ coincides
on $ F $ with $g_n$, and  let $x_{n+1}$ be a non-negative continuous
function supported on the subset of $V$
where $g_n$ attains its supremum on $ V$ up to
$\varepsilon$, i.e., on the set $\{t\in V\dopu g_{n}(t) > \sup_{s\in
V} g_{n}(s) - \eps\}$, etc. (There is no initial restriction in the
choice of $y_{1}$ and $x_{1}$  apart from the support and positivity
conditions.)

We first claim that
$$
\|g_{n}\|_{F} := \sup_{t\in F} |g_{n}(t)| \le 3+n\eps.
$$
This is certainly true for $n=1$ since $\|f_{1}+y_{1}\|< 1+\eps$.
Now induction yields for $t\in F$
\begin{align}
|g_{n+1}(t)|
&=
|g_{n}(t) + f_{n+1}(t)|
\nonumber\\
&=
\bigl| \|g_{n}\|_{F} \, y_{n+1}(t) + f_{n+1}(t) \bigr|
\nonumber\displaybreak[0]\\
&=
\bigl| y_{n+1}(t) + f_{n+1}(t) + ( \|g_{n}\|_{F} -1 ) y_{n+1}(t) \bigr|
\nonumber\displaybreak[0]\\
&\le
|y_{n+1}(t) + f_{n+1}(t)| + \|g_{n}\|_{F} -1
\nonumber\\
&\le
1+\eps + 2+n\eps = 3+ (n+1)\eps.   \nonumber
\end{align}
(We have tacitly assumed that $\|g_{n}\|_{F}\ge1$ since the induction
step  is clear otherwise, because $\|f_{n+1}\|\le 2+\eps$.)

Next, we  have that
$$
\sup_{t\in V}  g_{n}(t) > n(1-2\eps).
$$
Indeed, the functions $x_{1}$ and $y_{1}$ are disjointly supported;
hence by the definition
of $D(x_{1},y_{1}, \varepsilon)$ there is a point in
the support of $x_{1}$ at which  $f_{1} = g_{1}$ is
bigger than  $1 - \varepsilon$.
Thus, the above inequality holds for $n=1$.
To perform the induction step we use the same argument to find a
point $t_{0}$  in the support of $x_{n+1}$ at which $f_{n+1}$ exceeds
$1-\eps$.
At this point
$t_{0}$  the function $g_n$ attains its supremum on $V$ up to  $\varepsilon$. So
\begin{align}
\sup_{t\in V} g_{n+1}(t)
&\ge g_{n+1}(t_{0})
=
g_{n}(t_{0}) + f_{n+1}(t_{0})
>
\sup_{t\in V} g_n(t)  +  1 - 2\varepsilon
\nonumber \\
&> n(1-2\eps) + 1-2\eps
= (n+1)(1 - 2\varepsilon).  \nonumber
\end{align}

Therefore $\|g_{n+1} \| > (n+1)(1 - 2\varepsilon)$, and
on the other hand we have
$\|Tg_n\| \le \sum_{k=1}^n \|Tf_{n}\| <  n\varepsilon$.
So for $n$ big enough the function
$f=g_n/ \|g_n \|$ will satisfy the desired conditions.
\end{proof}

Actually, a somewhat smaller class of operators turns out to be
crucial.

\begin{definition} \label{defi3.5}
Let $X$ be a space with the Daugavet property.
Define the class of {\em narrow operators}\/
by $\narr(X) = \SDt(X) \mytilde - X^*$.
\end{definition}

In other words, an  operator $T$ is said to be a narrow operator
if, for every
$x^* \in X^*$, $T  \mytilde + x^*$ is a strong Daugavet operator.

Incidentally, it follows from the defining property of a narrow
operator that a Banach space on which at least one narrow operator is
defined must automatically have the Daugavet property.

Proposition~\ref{theo3.3} and
Proposition~\ref{prop2.7} imply that $\narr(X)$ is a $\sim$closed
homogeneous order ideal, and hence so is $\cp(\narr(X))$.

We now show that we won't get anything new on $C(K)$ if $K$ is
perfect.

\begin{theo} \label{cor3.6}
For a perfect compact Hausdorff space $K$,
the classes of $C$-narrow operators and of narrow
operators coincide on $C(K)$.
\end{theo}

\begin{proof}
Since a narrow operator is a strong Daugavet operator and a strong
Daugavet operator is $C$-narrow (Theorem~\ref{theo3.4}), it is left
to prove that a $C$-narrow operator $T$ on $C(K)$ is narrow if $K$ is
perfect. Let $x^*\in C(K)^*$ be a functional,
represented by a regular Borel
measure $\mu$; we have to show that $T \mytilde+ x^*$ is a strong
Daugavet operator.

Thus, let $f,g\in S(C(K))$ and $\eps>0$. Let $\eps'=\eps/(4+\|T\|)$,
and consider the open set $U=\{t\dopu |f(t)|>1-\eps'\}$. Pick an open
non-empty subset $V\subset U$ with the property that $f-g$ is almost
constant on $V$ in that for some number $c\in[-2,2]$
$$
|f(t)-g(t)-c|\le\eps' \text{ for }t\in V,
$$
and $|\mu|(V)\le\eps'$; the latter is possible since $K$  has no
isolated points. Since $T$ is $C$-narrow, there is some $\varphi\in
S(C(K))$ vanishing off $V$ such that $\|T\varphi\|\le\eps'$; in fact,
$\varphi$ can (and will) be chosen positive \cite[Lemma~1.4]{KadPop}.
Let $h= \varphi f + (1-\varphi) g$. Then $\|h\|\le1$, and 
$$
\|f+h\| \ge \sup_{t\in V} |f(t)+h(t)| \ge 2-2\eps' \ge 2-\eps;
$$
furthermore
\begin{align*}
|x^*(g) - x^*(h)| & =
|x^*(\varphi(f-g))| \le |\mu(V)| \,\|\varphi\| \, \|f-g\| \le 2\eps'
\\
\|T(g) - T(h)\| & =
\|T(\varphi(f-g))\| \le \|T\| \, \|\varphi(f-g-c)\| + \|T\varphi\| \,
|c| \\
& \le \|T\|\eps' + 2\eps' 
\end{align*}
so that 
$$
\| (T \mytilde+ x^*)(g-h) \| \le \eps,
$$
which proves that $T\mytilde+ x^*$ is a strong Daugavet operator.
\end{proof}

We shall show below (Section~\ref{sec6}) that in general, 
on a space with the Daugavet property
narrow and strong Daugavet operators are not the same.

We don't know if in general $\narr(X)$ is a subsemigroup
of ${\Op}(X)$ as it will be shown to be the case for $X=C(K)$ 
(Theorem~\ref{theo4.8}), but we will show that its central part
$\cp(\narr(X))$ is always large. It contains, in particular,
all strong Radon-Nikod\'ym  operators and all  operators which
do not fix copies of $\ell_1$. Hence all the operators
which are majorized by linear combinations of strong Radon-Nikod\'ym
operators and operators  not fixing copies of $\ell_1$,
as well as $\sim$-limits of sequences of such operators
belong to $\cp(\narr(X))$.

We now formulate a number of lemmas. 
Eventually, Proposition~\ref{Prop3.11} will present a
geometric description of what distinguishes a narrow operator from a
strong Daugavet operator.

\begin{lemma}\label{Lnarrow}
Let $T\in \narr(X)$. 
Then for every $x, y  \in S(X)$,
$\varepsilon > 0$  and every slice $S = S(x^*,\alpha)$
of the unit ball  of $X$ containing $y$ there is an element $v \in S$
such that  $\|x+v\|>2-\varepsilon$ and
$\|T(y-v)\|< \varepsilon$.
\end{lemma}

\begin{proof}
Fix some $0<\delta <\varepsilon$
and find an element $y_1$ in norm-interior of
$S$ such that $\|y-y_1\|<\delta$. By Proposition~\ref{prop1}, for every
$0<\delta_1 < \varepsilon$
there is an element
$u \in U_{x^*,\delta_1} \cap D(x ,y_1/\|y_1\|,\delta_1)$ such that
$\|Tu\|< \delta_1$. If $\delta_1$ is small enough, then
$v:= ( y_1+\|y_1\|u )/ \| (y_1+\|y_1\|u) \| \in S$.
So, if in turn $\delta$ is small enough, then $v $
satisfies our requirements.
\end{proof}

\begin{lemma}
For every  $\tau > 0$
and  every pair of positive numbers $a, b$ there is a
$\delta > 0$ such that if $v,x \in S(X)$ and
$\|x+v\|> 2 - \delta$,
then $\|ax+bv\|>a + b -\tau$.
\end{lemma}

\begin{proof}
Select $\delta <\tau / \max(a,b)$.
There are two symmetric cases: $a \le  b$ or $b \le  a$.
Consider for example the first of them. If we assume that our
statement is not true, then we obtain
\begin{align}
 2 - \delta &< \|x+v\|= \|(1-a/b)x + 1/b(ax+bv)\| \nonumber \\
            &\le 1-a/b + 1/b (a + b -\tau)= 2 - \tau/b,
\nonumber
\end{align}
a contradiction.
\end{proof}

\begin{lemma}\label{Lccs}
Let $T\in \narr(X)$. 
\begin{statements}
\item
Let $S_1, \dotsc ,S_{n}$ be a finite collection of slices
and $U \subset B(X)$ be
a convex combination of these slices, i.e., there are
$\lambda _k \ge 0$, $ k=1,\dotsc ,n$, $ \sum_{k=1}^n \lambda _k = 1$,
such that
$ \lambda _1S_1 + \cdots + \lambda _nS_n = U$.
Then for every  $\varepsilon > 0$, every $x_1 \in S(X)$
and  every  $w\in U$ there
exists an element $u \in U$ such that $\|u+x_1\|> 2 - \varepsilon$
and $\|T(w-u)\|<\varepsilon$.
\item
The same conclusion is true if $U$ is a relatively weakly open set.
\end{statements}
\end{lemma}

\begin{proof}
(a) First of all let us fix elements $y_k\in S_k$  such that
$ \lambda _1y_1 + \cdots + \lambda _ny_n = w$.
Applying repeatedly Lemma~\ref{Lnarrow} with sufficiently small
$\varepsilon_j$ to $S_j$, $y_j \in S_j$ and
$$
x_j = \biggl( x_1+ \sum_{k=1}^{j-1} \lambda _k v_k \biggr) \biggm/
\biggl\|x_1+ \sum_{k=1}^{j-1} \lambda _k v_k\biggr\|,
$$
we may select elements $v_k \in S_k$ with
$\|T(y_k-v_k)\|<\varepsilon$, $ k=1,\dotsc ,n$, in
such a way that for every $ j=1, \dotsc ,n$
$$
\biggl\| x_1+ \sum_{k=1}^j \lambda _k v_k \biggr\| >
1 + \sum_{k=1}^j \lambda _k(1 -\varepsilon)
$$
(to get the last inequality, we need to apply the
previous lemma at each step). The element
$u= \lambda _1v_1 + \lambda _2v_2+  \cdots + \lambda _n v_n$ will be
as required.

(b) This follows from (a) since given $u\in U$  there is a convex
combination $V$ of slices such that $u\in V\subset U$; see 
\cite[Lemma~II.1]{GGMS} or \cite{Shv1}.
\end{proof}

\begin{prop}\label{Prop3.11}
An  operator $T$ on a Banach space $X$ with the Daugavet property
is narrow if and only if  for every $x, y  \in S(X)$,
$\varepsilon > 0$  and every slice $S $
of the unit ball  of $X$ containing $y$ there is an element $v \in S$
such that  $\|x+v\|>2-\varepsilon$ and
$\|T(y-v)\|< \varepsilon$.
\end{prop}

\begin{proof}
It only remains to show that the above condition is sufficient for $T$
to be narrow. We first note that an operator satisfying that
condition will also satisfy the conclusion of Lemma~\ref{Lccs}; see
the proof of that lemma. Now, if $x_0^*\in X^*$, $x,y\in S(X)$ and
$\eps>0$, consider the relatively weakly open set
$$
U:= \{z\in B(X)\dopu |x_{0}^*(z-y)| < \eps/2 \}.
$$
By Lemma~\ref{Lccs} there exists some $w\in U$ such that
$\| w+x \| > 2-\eps/2$ and $\| T(w-y) \| < \eps/2$; note that $y\in
U$. By definition this means that $T \mytilde+ x_{0}^*$ is a strong
Daugavet operator; i.e., $T$ is narrow.
\end{proof}

Let $T$ be a strong Radon-Nikod\'ym  operator on  a space $X$ with
the Daugavet property; this means that the closure of $T(B(X))$ is a
set with the Radon-Nikod\'ym property, cf.\  \cite{Bou} for this notion.
We shall show that such an operator is narrow.
For $\varepsilon > 0$, consider the
subset $A(T,\varepsilon)$ of $B(X)$ defined by $y \in A(T,\varepsilon)$
if  there exists a convex combination $U$ of slices of
the unit ball such that $y \in U$ and $U \subset y + U_{T,\varepsilon}$.

\begin{lemma}
The set $A(T,\varepsilon)$ introduced above is a convex dense subset
of $B(X)$.
\end{lemma}

\begin{proof}
The convexity is evident. To prove the density we need to show, by the
Hahn-Banach theorem, that for every $x^* \in S(X^*)$ and
every $0<\delta <\varepsilon$ there is an element
$y \in A(T,\eps)$ such that $x^*(y) > 1 - \delta$ (in other words,
$y \in S=S(x^*,\delta)$).
Let us fix an element  $x \in B(X)$ with $x^*(x) > 1 - \delta/2$ and
consider  the operator $T_1= x^* \mytilde + T$. Consider further
$\overline{T_1(B(X))}$ and a $\delta/2 $-neighbourhood $W$ of $T_1x$ in
$\overline{T_1(B(X))}$. By the  Radon-Nikod\'ym property of the set
$\overline{T_1(B(X))}$
there is a convex combination $W_1$ of slices of $\overline{T_1(B(X))}$
in  $W$. The preimages in $B(X)$ of these
slices of $\overline{T_1(B(X))}$
are slices in $B(X)$. The corresponding convex
combination $U$ of these slices in $B(X)$ lies in the preimage
of $W$ in $B(X)$, so this convex combination is contained in
$(x + U_{T_1,\delta/2}) \cap  B(X)$. Fix an element $y \in U$.
By our construction
$y \in U \subset (x + U_{T_1,\delta/2}) \cap B(X)\subset S$.
On the other hand,
$$
U \subset x + U_{T_1,\delta/2}\subset y + U_{T_1,\delta}\subset
y + U_{T,\delta} \subset y + U_{T,\varepsilon},
$$
so $y \in A(T,\varepsilon)$.
\end{proof}

The following result is a generalisation of \cite[Th.~2.3]{KadSSW}.
It can be
understood as a transfer theorem: in Definition~\ref{defi3.5} one can
pass from one-dimensional operators to a much wider class of
operators. Let us denote the class of strong Radon-Nikod\'ym
operators on $X$ by $\SRN(X)$.

\begin{thm} \label{theo3.10}
Let $X$ be a space with the Daugavet property, $T$ be  narrow and
$T_{1}$ be a strong Radon-Nikod\'ym  operator on $X$. 
Then $T \mytilde + T_{1}$ is
narrow; that is, we have
$\narr(X) \mytilde+ \SRN(X) = \narr(X)$.
In particular every strong Radon-Nikod\'ym  operator $T_{1}$ on $X$ is a
narrow operator.
\end{thm}

\begin{proof}
 Let us fix $\varepsilon > 0$, $x,y \in S(X)$ and
$y_1 \in A(T_{1},\varepsilon)$ such that $\|y-y_1\|<\varepsilon$.
According to the definition of $A(T_{1},\varepsilon)$
there exists a convex combination $U$ of slices of
the unit ball such that $y_1 \in U$ and
$U \subset y_1 + U_{T_{1},\varepsilon}$. By Lemma~\ref{Lccs}
there is an element $z \in U$ such that
$\|z+x\|> 2 - \varepsilon$ and $\|T(y_1-z)\|<\varepsilon$.
But the inclusion $z \in y_1 + U_{T_{1},\varepsilon}$ means that
$\|T_{1}(y_1-z)\|<\varepsilon$. So
$$
\|(T\mytilde + T_{1})(y-z)\|<    \varepsilon \|T\mytilde + T_{1}\| +
\|(T\mytilde + T_{1})(y_1-z)\|<    \varepsilon \|T\mytilde + T_{1}\| +
2 \varepsilon.
$$
Because  $\varepsilon$ is arbitrarily small,
the last inequality shows that $T \mytilde + T_{1}$ satisfies
the definition of a strong Daugavet operator.

Now let $x^*\in X^*$ and consider $T_{2}= T_{1} \mytilde+ x^*$. This
is a strong Radon-Nikod\'ym operator, too. So $(T \mytilde+ T_{1})
\mytilde+ x^* = T \mytilde + T_{2}$ is a strong Daugavet operator by
what we have just proved; by definition, this says that $T\mytilde+
T_{1}$ is narrow.
\end{proof}

\begin{cor}\label{cpart}
Let $X$ be a Banach space with the Daugavet property.
\begin{statements}
\item
$\narr(X) \mytilde + X^* = \narr(X)$.
\item
$\cp(\narr(X)) =  \SDt(X) \mytilde - \narr(X)$.
\item
$\SRN(X) \subset \cp(\narr(X))$.
\end{statements}
\end{cor}

\begin{proof}
(a)
follows from the previous theorem, because every
finite-rank operator is a  strong Radon-Nikod\'ym  operator.

For (b) use Theorem~\ref{cpart1} and note that
\begin{align*}
\SDt(X) \mytilde - \narr(X)
&= \SDt(X) \mytilde - (\narr(X) \mytilde + X^*) \\
&=
(\SDt(X) \mytilde - X^*) \mytilde - \narr(X) \\
&=
\narr(X) \mytilde - \narr(X).
\end{align*}

(c)
is a restatement of Theorem~\ref{theo3.10}.
\end{proof}


\section{Operators which do not fix copies of $\ell_1$} \label{sec4}

It is proved in \cite{Shv1} that an operator $T\dopu X \to X$ on a
space with the Daugavet property which does not fix a copy of
$\ell_{1}$ satisfies the Daugavet equation. Recall that $T\in\Op(X)$ does not
fix   a copy of $\ell_{1}$ if there is no subspace $E\subset X$
isomorphic to $\ell_{1}$ on which the restriction $T\dopu E\to T(E)$
is an isomorphism.  By Rosenthal's $\ell_{1}$-theorem, this is
equivalent to saying that for every bounded sequence $(x_{n})\subset
X$, the sequence of images $(Tx_{n})$ admits a weak Cauchy
subsequence. We shall investigate the class of operators not fixing a
copy of $\ell_{1}$ in the present context.

We will use the following theorem, due to H.P.~Rosenthal
\cite{Ros-rec}:

\begin{thm}\label{ROS}
Let $X$ be a separable Banach space without $\ell_1$-subspaces.
If $A \subset X$ is bounded and $x^{**}\in X^{**}$ is a weak$\,^*$
limit point of $A$, then there is a sequence in  $A$ which
converges to $x^{**}$ in the weak$\,^*$ topology of $X^{**}$.
\end{thm}

In fact, we shall need a generalization of this result and first
provide a lemma.

\begin{lemma}\label{ROSlem}
Let $X$ be a Banach space without subspaces
isomorphic to $\ell_1$,  and let
$\{x_{n,m}\}_{n,m \in \N } \subset X$ be a bounded double sequence.
Let $x^{**}\in X^{**}$ be a $\sigma (X^{**},X^{*})$-limit point
of every column $\{x_{n,m}\}_{n \in \N }$ of
$\{x_{n,m}\}_{n,m \in \N }$. Then there are strictly increasing
sequences $(n(k))$, $( m(k))$ of indices such that
$x_{n(k), m(k)} \rightarrow x^{**}$ in $\sigma (X^{**},X^{*})$.
\end{lemma}

\begin{proof}
Consider an auxiliary space $Y=X \times \R$ and an auxiliary
matrix $\{y_{n,m}\}_{n,m \in \N } \subset Y$,
$y_{n,m}=(x_{n,m}, 1/n + 1/m)$. Since $Y$ contains
no copies of $\ell_1$ either and since $(x^{**}, 0)$ is a $\sigma
(Y^{**},Y^{*})$-limit point of $\{y_{n,m}\}_{n,m \in \N }$,
there is, according to Theorem~\ref{ROS}, a sequence of the
form $(y_{n(k), m(k)})$
which converges to $(x^{**}, 0)$ in $\sigma (Y^{**},Y^{*})$.
This means in particular that
$x_{n(k), m(k)} \rightarrow x^{**}$ in $\sigma (X^{**},X^{*})$
and $1/n + 1/m \rightarrow 0$. So $(n(k))$ and $( m(k))$ both tend to
$\infty$, which, after passing to a
subsequence, provides the desired sequence.
\end{proof}

The next result is a direct generalisation  of
Theorem~\ref{ROS}.

\begin{thm}\label{ROS1}
Let $X$ be a separable Banach space without $\ell_1$-subspaces,
$(\Gamma, {\preceq})$ be a directed set, and let
$F\dopu \Gamma \rightarrow X$
be a bounded function. Then for every $\sigma (X^{**},X^{*})$-%
limit point $x^{**}$ of the function $F$  there is a strictly
increasing sequence
$\gamma(1) \preceq \gamma(2) \preceq \ldots$
in $\Gamma$ such that $\bigl( F(\gamma(n)) \bigr)$ converges to  $x^{**}$ in
$\sigma (X^{**},X^{*})$.
\end{thm}

\begin{proof}
Using inductively
Theorem~\ref{ROS} we can select a doubly indexed sequence
$\{\gamma_{n, m}\}_{n,m \in \N }$
in $\Gamma$ with the following properties:
\begin{enumerate}
\item
for  every $m \in \N$, $x^{**}\in X^{**}$ is a
$\sigma (X^{**},X^{*})$-limit point
of every column $\{F(\gamma_{n,m})\}_{n \in \N }$;
\item
for  every $m, n, k, l  \in \N$, if $\max\{k, l\} < m$,
then $\gamma_{k,l} \preceq \gamma_{n,m}$.
\end{enumerate}

Applying Lemma~\ref{ROSlem} and passing to a subsequence if necessary,
we obtain strictly increasing
sequences $(n(k))$, $( m(k))$ such that
$\max_{k<j} \{n(k), m(k)\} < m(j)$ and $ \bigl( F(\gamma_{n(k), m(k)}
\bigr)$
converges to  $x^{**}$ in $\sigma (X^{**},X^{*})$.
To finish the proof put
$\gamma(k) = \gamma_{n(k), m(k)}$.
\end{proof}

We wish to prove that an operator not fixing a copy of $\ell_{1}$ is
narrow (Theorem~\ref{l1}  below). To cover the case of non-separable
spaces
as well we first show that the Daugavet property is separably
determined. The next lemma prepares this result.

\begin{lemma} \label{finit}
Let $X$ be a Banach
space with the Daugavet property. Then for any $\eps>0$
and $x,y\in S(X)$, there exists a finite-dimensional subspace
$Y=Y(x,y,\eps)$ of $X$ with $x,y\in Y$
such that for every slice $S(x^*,\eps/2)$
containing $y$ there is some $y_{1}\in S(Y)\cap S(x^*,\eps)$ such
that $\|y_{1}+x\|>2-\eps$.
\end{lemma}

\begin{proof}
Assume there exist $\eps>0$ and $x,y\in S(X)$ such that for every
finite-dimensional subspace $Y\subset X$ there is a slice
$S(x_{Y}^*, \eps/2)$ containing $y$ with $\|y_{1}+x\|\le 2-\eps$ for
all $y_{1}\in S(Y)\cap S(x_{Y}^*,\eps)$. Take a weak$^*$ cluster
point $x^*$ of the net $(x_{Y}^*)$ and let $x_{0}^*= x^*/\|x^*\|$. We
have $x^*(y) \ge 1-\eps/2$ since $y\in S(x_{Y}^*,\eps/2)$ and
therefore $\|x^*\|\ge 1-\eps/2$. Now if $y_{1}\in S(x_{0}^*,\eps/2)$,
then $x^*(y_{1})\ge \|x^*\| (1-\eps/2) >1-\eps$ and therefore
$x_{Y_{1}}^*(y_{1})>1-\eps$ for some $Y_{1}$ that contains $y_{1}$.
So by assumption $\|y_{1}+x\|\le 2-\eps$, which contradicts
Lemma~\ref{lem1.1} when applied to the slice $S(x_{0}^*, \eps/2)$.
\end{proof}

\begin{thm}\label{separable}
A Banach space $X$ has the Daugavet property if and only if
for every separable subspace $Y \subset X$
there is a  separable subspace $Z \subset X$ which
contains $Y$ and has the Daugavet property.
\end{thm}

\begin{proof}
Suppose $X$ has the Daugavet property.
Let $(v_n)$ be a dense sequence in $Y$. We select a sequence
$V_1 \subset V_2 \subset \ldots$
of finite-dimensional subspaces of $X$ by the following inductive
procedure. Put  $V_1 = \Lin v_1$. Suppose $V_n$ has already been constructed.
Fix a $2^{-n}$-net $(x_k^n, y_k^n)$, $ k=1,  \dotsc, N_n$, in
$S(V_n) \times S(V_n)$ provided with the sum norm,
select by Lemma~\ref{finit}
finite-dimensional subspaces
$Y_k=Y(x_k^n,y_k^n,\eps)$, $ k=1,  \dotsc, N_n$,
for $\eps = 2^{-n}$ and define
$V_{n+1}= \Lin ( \{v_{n+1}\}\cup Y_1 \cup \dotsc \cup Y_{N_n} )$.

If $Z$ is defined to be the closure of the union of all the $V_{n}$,
then $Y\subset Z$ and $Z$ has the Daugavet property by
Lemma~\ref{lem1.1}.

Conversely, let $x\in S(X)$, $\eps>0$ and let $S\subset B(X)$ be a
slice. Fix a point $z\in S$. If $Z$ is a separable subspace with the
Daugavet property containing $x$ and $z$, then by Lemma~\ref{lem1.1}
there exists some $y\in S\cap Z$  such that $\|y+x\|>2-\eps$. Again
by Lemma~\ref{lem1.1} this shows that $X$ has the Daugavet property.
\end{proof}

We shall need the operator version of this theorem, which is based on
the following lemma. The proofs of Lemma~\ref{finit2} and
Theorem~\ref{sep} are virtually the same as those
of Lemma~\ref{finit} and Theorem~\ref{separable} (one uses
Proposition~\ref{Prop3.11}).

\begin{lemma} \label{finit2}
Let $X$ be a Banach
space with the Daugavet property and let $T$ be a narrow
operator on~$X$. Then for any $\eps>0$
and $x,y\in S(X)$, there exists a finite-dimensional subspace
$Y=Y(x,y,\eps)$ of $X$ with $x,y\in Y$
such that for every slice $S(x^*,\eps/2)$
containing $y$ there is some $y_{1}\in S(Y)\cap S(x^*,\eps)$
with $\|Ty_{1}-Ty\| < \eps$ such
that $\|y_{1}+x\|>2-\eps$.
\end{lemma}

\begin{thm}\label{sep}
An operator $T$ on a Banach space $X$ is narrow if and only if for
every separable subspace $Y$ of $X$ there is a separable subspace
$Z\subset X$ containing $Y$ such that the restriction of $T$  to $Z$
is narrow.
\end{thm}

This theorem leads to an important structural result on narrow
operators on $C(K)$.

\begin{theo}  \label{theo4.8}
If $K$ is a perfect compact Hausdorff space, then $\narr(C(K))$ is a
subsemigroup of $\Op(C(K))$.
\end{theo}

\begin{proof}
First, let $K$ be a perfect compact metric space.
It follows from Theorem~\ref{cor3.6} and \cite[Th.~1.8]{KadPop} that
the  set of all  narrow
operators on $C(K)$
is stable under the operation $ \mytilde +$,
i.e., it is a semigroup. (In fact, \cite{KadPop} only deals with
$K=[0,1]$, but the arguments work as well for a metric~$K$.)

We shall now reduce the general case to the metric one.
Let now $K$ be a perfect compact Hausdorff space, and
let $T_{1}$ and $T_{2}$ be two narrow operators on $C(K)$;
we shall verify that $T_{1} \mytilde+ T_{2}$
is narrow, using Theorem~\ref{sep} above.

Thus, let $Y$   be a separable subspace of $C(K)$.
We shall first argue that there is a separable space $Z_{1}$
containing $Y$ such that $\rest{T_{1}}{Z_{1}}$ and
$\rest{T_{2}}{Z_{1}}$ are strong Daugavet  operators. 
Let $A$ be a countable dense subset of $S(Y)$. For every pair $(x,y)$
in $A\times A$ and every $\eps= 1/k$ there is some $z_{1}$
(resp.~$z_{2}$) according to the definition of the strong Daugavet 
property of $T_{1}$ (resp.~$T_{2}$). The countable collection of
these $z$'s and $Y$ span a closed separable subspace $X_{1}$.
Repeating this procedure starting from $X_{1}$ yields some
closed separable subspace $X_{2}$, etc. 
The closed linear span $Z_{1}$ of $X_{1}, X_{2}, X_{3}, \ldots$
then has the desired property.

Now by \cite[Lemma~2.4]{WeisDirk} there is a separable space
$Z_{2}\supset Z_{1}$ isometric to some space $C(M_{2})$ for a
perfect compact metric space~$M_{2}$. By the same token as above, we
can extend $Z_{2}$ to a separable space $Z_{3}$ so that $T_{1}$ and
$T_{2}$ are strong Daugavet  operators on $Z_{3}$, and we 
can extend $Z_{3}$ to a separable space $Z_{4}$ 
isometric to some space $C(M_{4})$      for a
perfect compact metric space~$M_{4}$, etc.
Let $Z$ be the closed linear span of $Z_{1}, Z_{2}, Z_{3}, \dotsc$.
Then $\rest{T_{1}}{Z}$  and $\rest{T_{2}}{Z}$ are strong Daugavet 
operators, and $Z$ is 
isometric to some space $C(M)$  for a
perfect compact metric space~$M$. By what we already know, 
$\rest{(T_{1} \mytilde+ T_{2})}{Z}$
is a narrow operator on $Z\cong C(M)$; recall that the classes of
narrow and strong Daugavet operators coincide on $C(M)$. 

Finally, Theorem~\ref{sep} implies that $T_{1} \mytilde+ T_{2}$
is narrow on $C(K)$, which proves the theorem.
\end{proof}

Next we introduce a topology related to an order ideal of operators.

\begin{definition}
Let ${\cal M}\subset\Op(X)$ be an order ideal of operators, closed
under the operation $ \mytilde +$. Then the system of tubes
$U_{T,\varepsilon}$, $ T \in {\cal M}$, $\varepsilon > 0$, defines a base
of neighbourhoods of\/ $0$ for some locally convex topology on $X$.
We denote this topology by $\sigma (X, {\cal M})$.
\end{definition}

If ${\cal M}={\cal F}(X)$,  the class of all finite-rank operators, then
$\sigma (X, {\cal M})$ coincides with the weak topology; if
${\cal M}=\Op(X)$, then
$\sigma (X, {\cal M})$ coincides with the norm topology. For
classes which are in between one gets  topologies which are
 between the weak and the norm topology. If  ${\cal N}$ is
a  collection of subsets in $X$ such that ${\cal N}^{\sim}$ is  closed
under the operation $ \mytilde +$, then $\sigma (X, {\cal  N}^{\sim})$
is the strongest  locally convex
topology on $X$ continuous with respect to the norm, in which the zero vector
belongs to the closure of every element of  ${\cal N}$.

\begin{definition}
A locally convex topology $\tau$ on $X$
is said to be a  {\em Daugavet topology}\/ if for every two elements
$x, y \in S(X)$, for every
$\varepsilon > 0$ and every $\tau$-neighbourhood $U$ of $y$ there is
an element
$z \in  U \cap  S(X)$ such that $\|z+x\|> 2 - \varepsilon$.
\end{definition}

Of course, $\sigma (X, {\cal M})$ is  a  Daugavet topology if and only
if every operator $T \in {\cal M}$ is a strong Daugavet operator.

\begin{lemma}
Let $X$ be a Banach space with the Daugavet property, $T$
a narrow operator, $A=\{a_1, \dots , a_n\} \subset S(X)$,
$\varepsilon > 0$ and  $y \in S(X)$. Then for every
$\sigma (X,\cp(\narr(X)))$-%
neighbourhood $W$ of $y$ there is an element $w \in W \cap  S(X)$
such that $\|T(w-y)\|<\varepsilon$ and
$\|w+a\|> 2 - \varepsilon$ for every $a \in A$.
\end{lemma}

\begin{proof}
We shall argue by induction on $n$. First of all consider
$n=1$. Every $\sigma (X,\cp(\narr(X)))$-neighbourhood
 of $y$ can be represented
as $W=y + U_{R,\delta}$, where $R \in \cp(\narr(X))$.
Since $T_1 = R \mytilde + T$ is a strong Daugavet operator by
definition of the central part,
there is an element $w \in S(X)$ such that
$\|w+a_1\|> 2 - \varepsilon$ and
$\|T_1(w-y)\|< \min(\delta, \varepsilon)$.
The last inequality means, in particular, that
$\|T(w-y)\|<\varepsilon$ and $w \in W$.

Now suppose our assertion is true for $n$, let us
prove it for $n+1$. Let
$A=\{a_1, \dots , a_n, a_{n+1}\} \subset S(X)$,
and let us assume that an element $w_1 \in W \cap  S(X)$
such that $\|T(w_1-y)\|<\varepsilon/2$ and
$\|w_1+a_k\|> 2 - \varepsilon$, $k=1, \dots, n$,
has already been selected. Then there is a weak  neighbourhood
$U$ of $w_1$ such that the inequalities
$\|u+a_k\|> 2 - \varepsilon$, $k=1, \dots, n$,
hold for every $u \in U$. The intersection
$U \cap W$ is  a  $\sigma (X,\cp(\narr(X)))$-neighbourhood
of $w_1$, so according to our inductive assumption for $n=1$,
there is an element $w \in S(X) \cap U \cap W$ such that
$\|w+a_{n+1}\|> 2 - \varepsilon$ and $\|T(w-w_1)\|<\varepsilon/2$.
This element $w$ satisfies all the requirements.
\end{proof}

Using an $\varepsilon$-net of the unit ball of the finite-dimensional
subspace $Z$ below one can easily deduce the following corollary.

\begin{prop}\label{weakl1}
Let $X$ be a Banach space with the Daugavet property, $T$ be
a narrow operator and $Z \subset X$ be a
finite-dimensional subspace. Then for every  $\varepsilon > 0$,
every  $y \in S(X)$ and every
$\sigma (X,\cp(\narr(X)))$-neighbourhood $W$ of $y$
there is an element $w \in W \cap  S(X)$
such that $\|T(w-y)\|<\varepsilon$ and
$\|z+w\|> (1 - \varepsilon)(\|z\|+\|w\|)$ for every
$z \in Z$.
\end{prop}


\begin{thm}\label{l1}
Let $X$ be a Banach space with the Daugavet property and let $T$ be
an operator on $X$ which does not fix a copy of  $\ell_1$.
Then $T \in \cp(\narr(X))$, so in particular $T$ is a narrow operator.
\end{thm}

\begin{proof}
Lemma~1(xii) of \cite{DFJP} implies that
every  operator which does not fix a copy of  $\ell_1$
can be factored through a space without $\ell_1$-subspaces. So
every  operator which does not fix a copy of  $\ell_1$ can be majorized
by an operator which maps into  a space without $\ell_1$-subspaces.
Since the class of narrow operators is an order ideal, it is enough
to prove our theorem for $T\dopu  X \rightarrow Y$, where $Y$ has no
 $\ell_1$-subspaces. Also, by Theorem~\ref{sep} we may assume
that  $X$ and $Y$ are  separable.

Let us fix a narrow operator $R$, $\varepsilon > 0$  and
$x,y \in S(X)$. Let us introduce a directed  set $(\Gamma, {\preceq})$
as follows: the elements of  $\Gamma$ are finite sequences in $S(X)$ of the
form $\gamma = (x_1, \dotsc, x_n)$, $n \in \N$, with $x_1 = x$.
The (strict) ordering is defined by
$$
(x_{1},\dotsc,x_{n}) \prec (y_{1},\dotsc,y_{m}) \iff
n< m \ \& \ \{x_{1},\dotsc,x_n\} \subset \{y_{1},\dotsc,y_{m-1}\}
$$
and of course $\gamma_1 \preceq \gamma_2$ if 
$\gamma_1 \prec \gamma_2$ or $\gamma_1 = \gamma_2$.
Now define
a bounded function $F\dopu \Gamma \rightarrow Y \times \R \times \R$
by
$$
F(\gamma) = (Tx_n, \alpha(\gamma), \|R(y - x_n)\|),
$$
where
$$
\alpha(\gamma)=
\sup \bigl\{ a>0\dopu  \|z + x_n \| > a(\|z\| + \|x_n \|)
\ \forall z \in \Lin \{x_1, x_2, \dots, x_{n-1}\}
\bigr\}.
$$
Due to Proposition~\ref{weakl1}, for every weak neighbourhood $U$ of
$y$ in  $B(X)$, every  $\varepsilon > 0$ and every finite collection
$\{v_1, \dotsc , v_n\} \subset X$ there is some $v_{n+1} \in U$ for which
$\alpha \bigl( (v_1, \dotsc, v_{n+1}) \bigr) > 1 - \varepsilon$
and $\|R(y - v_{n+1})\| < \varepsilon$. This means that
$(Ty, 1, 0)$ is a weak
limit point of the function $F$. So, by Theorem~\ref{ROS1}
there is a strictly $\prec$-increasing
sequence $(\gamma_j) = \bigl( (x_1, \dotsc,x_{n(j)})
\bigr)$ for which $(Tx_{n(j)})$ tends weakly to $Ty$,
$(\|R(y - x_{n(j)})\|)$ tends to 0 and $(\alpha(\gamma_j))$ tends
to~$1$. Passing to a subsequence we can select points
$x_{n(j)}$ in such a way
that the sequence $\{x, x_{n(1)}, x_{n(2)},  \dots\}$ is
$\varepsilon$-equivalent to the canonical basis of  $\ell_1$. According to
Mazur's theorem, there is a sequence $z_n \in \conv\{x_{n(j)}\}_{j>n}$
such that $\|Ty-Tz_n\| \rightarrow 0$. Evidently
$\|z_n+x\|> 2 - \varepsilon$ and
$\| (R \mytilde +T)(y-z_n)\| \rightarrow 0$,
which means that $R \mytilde +T \in \SDt(X)$ and thus proves the theorem
by Corollary~\ref{cpart}(b).
\end{proof}

There are other applications of Theorem~\ref{ROS1}
which are not related to the  Daugavet property.
As an example let us prove the following theorem which was
earlier established by E.~Behrends \cite{Beh-bar} under the more restrictive
condition of separability of $X^*$.

\begin{thm}
Let $X$ be a Banach space without $\ell_1$-subspaces and
$A_n \subset X$ be bounded subsets with $0\in \clconv A_{n}$ for each
$n\in\N$.
Then there exists a sequence $(a_n)$ in $ X$ with $ a_n \in A_n $
for every $n$ such that $0\in \clconv \{a_{1},a_{2},\dots\}$.
\end{thm}

\begin{proof}
In each $A_n$ there is a separable subset whose closed convex hull
contains $0$. So,
passing to the linear span of these separable
subsets we may assume that  $X$ is  separable.
Introduce a directed set $(\Gamma, {\preceq})$
as follows: the elements of  $\Gamma$ are of the
form
$$
\gamma = \bigl( n, m, \{a_k\}_{k=n}^m, \{\lambda_k\}_{k=n}^m \bigr),
$$
where $n, m \in \N$, $ n < m$, $ a_k \in A_k $, $\lambda_k > 0$,
$\sum_{k=n}^m \lambda_k = 1$. Define $ \preceq$ as follows:
let $\gamma_1 = (n_1, m_1, \{a_k\}_{k=n_1}^{m_1},
\{\lambda_k\}_{k=n_1}^{m_1})$,
 $\gamma_2 = (n_2, m_2, \{b_k\}_{k=n_2}^{m_2}, \{\mu_k\}_{k=n_2}^{m_2})$;
then $\gamma_1 \preceq \gamma_2$ if $ m_1 < n_2$.
Define $F\dopu \Gamma \rightarrow X$ by the formula
$F(\gamma) = \sum_{k=n}^m \lambda_k  a_k$.
Now, $0$ is a weak limit point of $F$; see the proof of
\cite[Th.~4.3]{Beh-bar}.
So, by Theorem~\ref{ROS1} there is a sequence of elements
$$
\gamma_j = (n_j, m_j, \{a_k\}_{k=n_j}^{m_j}, \{\lambda_k\}_{k=n_j}^{m_j})
$$
such that $n_1 <  m_1 < n_2 <  m_2 < n_3 < \ldots$
and $\sum_{k=n_j}^{m_j} \lambda_k  a_k$ tends weakly to zero.
To finish the proof one just needs to apply Mazur's theorem.
\end{proof}


\section{Rich subspaces} \label{sec5}

In \cite{PliPop} a subspace $Y$ of $L_{1}$ is called rich if the
quotient map $q\dopu L_{1}\to L_{1}/Y$ is $L_{1}$-narrow, and
likewise a subspace $Y$ of $C(K)$ is called rich in \cite{KadPop} if the
quotient map $q\dopu C(K)\to C(K)/Y$ is $C$-narrow. We are now in a
position to discuss rich subspaces in general.

\begin{definition}
Let $X$ be a Banach space with the Daugavet property. A subspace $Y$ is
said to be {\em almost rich}\/
if the quotient map $q\dopu  X \rightarrow X/Y$ is
a strong Daugavet operator. A subspace $Y$ is
said to be {\em rich}\/ if  the quotient map $q\dopu  X \rightarrow X/Y$ is
a narrow operator.
\end{definition}

By Theorem~\ref{cor3.6} the new definition comprises the old one
for subspaces of $C(K)$.

The necessity to distinguish rich and almost rich subspaces will
become apparent later when we show that the following theorem does
not extend to almost rich subspaces; see Theorem~\ref{theo6.7}.

\begin{thm} \label{theo5.3}
A  rich  subspace $Y$ of  a Banach  space  $X$ with the Daugavet
property has the
Daugavet  property  itself. Moreover, $(Y,X)$ is a Daugavet
pair.
\end{thm}

\begin{proof}
 Consider elements $x\in S(X)$, $y \in S(Y)$,
a slice $S=S(x^*,\varepsilon)$ and  $y \in S$. According
to our assumption the quotient map $q\dopu  X \rightarrow X/Y$ is
a narrow operator. So there is an element  $u \in S$ such that
$\|u+x\|> 2 - \varepsilon$ and $\|q(y-u)\|=\|q(u)\| <\varepsilon$.
The last condition means that the distance from $u$ to  $Y$
is smaller than $\varepsilon$, so there is an element  $v \in Y$
with  $\|v-u\| <\varepsilon$.
The norm of $v$ is close to $1$, viz.\
$1-2 \varepsilon < \|v\|  < 1 + \varepsilon$. Put $w=v/\|v\|$.
For this $w$ we have $\|w-u\| < 3\varepsilon$, so
$w \in S(x^*,4 \varepsilon)$ and $\|w+x\|> 2 - 4\varepsilon$.
\end{proof}

This theorem leads to new hereditary properties for the Daugavet
property.

\begin{prop} \label{cor5.5}
Suppose $Y$ is a subspace  of a Banach space $X$ with the Daugavet property.
\begin{statements}
\item
If the quotient space $X/Y$ has
the Radon-Nikod\'ym property, then $Y$ is rich.
\item
If the quotient space $X/Y$ contains no copy of $\ell_{1}$,
then $Y$ is rich in $X$.
\item
If $(X/Y)^*$ has
the Radon-Nikod\'ym property, then $Y$ is rich.
\end{statements}
In either case $Y$ has the Daugavet property itself.
\end{prop}

\begin{proof}
(a) follows from Theorem~\ref{theo3.10}, (b) from
Theorem~\ref{l1}, and (c) follows from~(b).
\end{proof}

That $Y$ has the Daugavet property under assumption (a)
has been proved earlier in \cite{Shv1}.

\begin{rem}
If  the quotient map $q\dopu  X \rightarrow X/Y$ belongs to
$\cp(\narr(X))$,
then the restriction to  $Y$  of every
narrow operator on $X$ is a narrow operator itself.
If $Y$ is a rich  subspace of  a   space  $X$
having the Daugavet property,
then the restriction to  $Y$  of every
operator $T \in \cp(\narr(X))$ is a narrow operator.
\end{rem}

\begin{definition} \label{defi5.6}
We say that a subspace  $Y$ of  a 
space $X$ with the Daugavet property is {\em wealthy}\/ 
if $Y$ and every  subspace of $X$ containing $Y$ have
the Daugavet property.
\end{definition}

It is plain that if $Y$ is an (almost) rich subspace of  a 
space $X$ with the Daugavet property,
 then every bigger subspace is  (almost) rich, too. Thus, if
$Y$ is rich, then it is wealthy. We now investigate
the converse implication.

\begin{lemma} \label{lem5.7}
The following conditions for a subspace  $Y$ of  a Banach
space $X$ with the Daugavet property are equivalent:
\begin{aequivalenz}
\item
$Y$ is wealthy.
\item
Every  finite-codimensional
subspace of  $Y$ is wealthy.
\item
For  every pair
$x, y \in S(X)$, the linear span of  $Y$, $ x$ and $y$ has
the Daugavet property.
\item
For  every $x, y \in S(X)$, for every  $\varepsilon > 0$
and for every slice $S$ of $ S(X)$
which contains $y$ there is an element
$v  \in \Lin(\{x,y\}\cup Y) \cap S$ such that $\|x+v\|>2-\varepsilon$.
\end{aequivalenz}
\end{lemma}

\begin{proof}
Due to Proposition~\ref{cor5.5} every finite-codimensional subspace of
a  space with the Daugavet property has the Daugavet property itself
(see also \cite[Th.~2.14]{KadSSW});
this is the reason for the equivalence
of (i) and~(ii).
The implication (i) $ \Rightarrow $ (iii) follows immediately from
the definition of  a wealthy  subspace; (iii) $ \Rightarrow $ (iv)
and (iv) $\Rightarrow $ (i)
are consequences of Lemma~\ref{lem1.1}.
\end{proof}

Let us say that a pair of elements  $x, y \in S(X)$ is
{\em $\varepsilon$-fine}\/ if there is a slice $S$ of $ S(X)$ which contains
$y$ and the diameter of $S \cap \Lin\{x,y\}$ is less
than $\varepsilon$.

\begin{lemma} \label{lemma5.9}
Let $Y$ be a wealthy subspace of  a Banach
space $X$ with the Daugavet property
 and let a pair  $x, y \in S(X)$ be $\varepsilon$-fine.
Then  $Y$ intersects   $D(x ,y, 2 \varepsilon)$.
\end{lemma}

\begin{proof}
First of all let us fix a slice $S=S({x^*, \varepsilon_1})$
from the definition
of an $\varepsilon$-fine pair and  fix a $\delta > 0$ such that the set
$W = \{w \in \Lin\{x,y\}\dopu  \|w \| < 1+\delta$,
$ x^*(w)> 1-\varepsilon_1 \}$ still has
diameter less  than $\varepsilon$. Now let us find a finite-%
codimensional subspace $E \subset Y$ such that
\begin{enumerate}
\item
$ x^*=0$ on  $E$,
\item
if $e\in E$ and $ w \in \Lin\{x,y\}$, then
$ \|w \| < (1+\delta) \|e + w  \|$;
\end{enumerate}
the last condition can be satisfied by a variant of the
Mazur argument leading to the basic sequence
selection principle; see \cite[Lemma~6.3.1]{KadKad}. According to our
assumptions $\Lin(\{x,y\}\cup E)$ has the  Daugavet property. So there is
an element $v \in \Lin(\{x,y\}\cup E) \cap S$ such that
$\|x+v\|>2-\varepsilon$. Let us represent $v$ in the form
$v=e + w$, where  $e\in E$, $ w \in \Lin\{x,y\}$. By choice of $E$
this means that $ \|w \| < 1+\delta$ and
$x^*(w)=x^*(v)>1-\varepsilon_1$.
Thus, $w \in W$ and $ \|y - w \|< \varepsilon$.
Finally we have that the element $e$ belongs to
$E \cap D(x ,y,2 \varepsilon)$, which concludes the proof.
\end{proof}

Let us recall the following result \cite{KadSSW}, which can
also be deduced from our Proposition~\ref{weakl1}:

\begin{lemma}\label{L2}
Let $X$ be a Banach
space with the Daugavet property and  $Z \subset X$ be a
finite-dimensional subspace. Then, for every  $\varepsilon > 0$ in
every slice of the unit sphere of $X$ there is an element $x$
such that
\[
\|z+x\|> (1 - \varepsilon)(\|z\|+\|x\|) \qquad\forall
z \in Z.
\]
\end{lemma}

We now present two easy lemmas.

\begin{lemma} \label{lem5.2}
A subspace $Y$ of a  Banach space with the Daugavet property
which
is almost rich together with all of its
$1$-codimensional subspaces is rich.
\end{lemma}

\begin{proof}
Let $q\dopu X\to X/Y$ be the quotient map and let $x^*\in S(X^*)$;
further let $Y_{1} = Y\cap \ker x^*$ and let $q_{1}\dopu X \to
X/Y_{1}$ be the corresponding quotient map. Then $Y_{1}=Y$ or $Y_{1}$
is $1$-codimensional in $Y$. Now, in either case  we have
$\|q(x)\|+|x^*(x)|\le 2\|q_{1}(x)\|$ for all $x\in X$.
Since $q_{1}$ is a strong Daugavet
operator by assumption, so is $q\mytilde+ x^*$, and $q$ is narrow.
\end{proof}

\begin{lemma}\label{L1}
A subspace $Y$ of a  Banach space $X $ 
with the Daugavet property is almost rich if and only if  $Y$
intersects all the elements of ${\cal D}(X)$.
\end{lemma}

\begin{proof}
If $Y$
intersects all the elements of ${\cal D}(X)$, then the
quotient map $q\dopu  X \rightarrow X/Y$ is unbounded from below
on every element of ${\cal D}(X)$. So the quotient map belongs
to ${\cal D}(X)^{\sim}$ which coincides with the class of strong
Daugavet operators by Proposition~\ref{theo3.3}.

Now consider the converse statement.
If $Y$ is almost rich, then for every $\varepsilon > 0$ the map
$q$ is unbounded from below on every set of the form
$ D(x ,y, \varepsilon /2)$. This means that there is an element
$z \in Y$ for which
$\dist(z, D(x ,y, \varepsilon /2)) <  \varepsilon /2$. In this case
$z$ belongs to  $ D(x ,y, \varepsilon)$, so the intersection of this
set with $Y$ is non-empty.
\end{proof}

The following is the key result for establishing that wealthy
subspaces are rich.

\begin{lemma} \label{theo5.10}
Every  wealthy subspace $Y$  of  a Banach 
space $X$ having the Daugavet property is almost rich.
\end{lemma}

\begin{proof}
According to Lemma~\ref{L1} we need to prove that for every  positive
$\varepsilon < 1/10$ and every pair  $x, y \in S(X)$ the subspace
$Y$ intersects   $D(x ,y, \varepsilon)$. To do this, according
to Lemma~\ref{lemma5.9}, it is enough to show that for every
$\varepsilon > 0$ and every pair  $x, y \in S(X)$ there is
an $\varepsilon$-fine pair  $x_1, y_1 \in S(X)$ which approximates
$(x, y)$ well; i.e., $\|x-x_1\|+\|y - y_1 \|< \varepsilon$.
Let us fix a positive $\delta <\varepsilon^2/8$ and select an
element $z \in S(X)$  in such
a way that for every $ w \in \Lin\{x,y\}$ and for every $t > 0$
$$
\|w+tz\| \ge(1-\delta)( \|w\|+|t|)
$$
(we use Lemma~\ref{L2}).
Put $x_1=x + \varepsilon z$, $ y_1=y$. To show that  $(x_1, y)$ is
an $\varepsilon$-fine pair it is sufficient to demonstrate that,
for every $v \in \Lin\{x_1,y\}$ with $\|v\| \ge \varepsilon$,
$\max\{\|y+v\|,\|y-v\| \}>1$. To do this let us argue ad absurdum.
Take some $v=ay+b(x+\varepsilon z)$ with $ \|v\| \ge \varepsilon$
and assume that $\max\{\|y+v\|,\|y-v\| \}=1$. Then
\begin{align*}
1&= \max\{\|y+ay+b(x+\varepsilon z)\|,\|y-ay-b(x+\varepsilon z)\| \} \\
&\ge
(1-\delta)(\max\{\|y+ay+bx+\|,\|y-ay-bx\| \}+|b|\varepsilon )\\
&\ge
(1-\delta)(1+|b|\varepsilon).
\end{align*}
So $|b|\le \varepsilon/4$. But in this case $|a|> \varepsilon/2$
and
$$
\max\{\|y+v\|,\|y-v\| \} > \max\{\|y+ay\|,\|y-ay\| \} - \varepsilon/3 >
1 + \varepsilon/6,
$$
which provides a contradiction.
\end{proof}

\begin{theo}
The following properties of a subspace  $Y$ of  a Daugavet
space $X$ are equivalent:
\begin{aequivalenz}
\item
$Y$ is wealthy.
\item
$Y$ is rich.
\item
Every finite-codimensional subspace of $Y$ is rich.
\end{aequivalenz}
\end{theo}

\begin{proof}
It is clear that (iii) $\Rightarrow$ (ii) $\Rightarrow$ (i), see the
remark following Definition~\ref{defi5.6}. Now suppose~(i).
Every $1$-codimensional subspace of $Y$ is wealthy by Lemma~\ref{lem5.7}
and is hence almost rich by Lemma~\ref{theo5.10}. An appeal to
Lemma~\ref{lem5.2} completes the proof.
\end{proof}


\section{Operators on $L_1$}  \label{sec6}

In this section we shall study strong Daugavet and narrow operators
on $L_{1}$. We first introduce a technical definition.

Let $(\Omega, \Sigma, \mu)$  be an atomless probability space.
A function $f \in L_1 = L_1(\mu)$ is said to be a
{\em   balanced $\eps$-peak}\/ on $A \in \Sigma$ if $f\ge-1$, $\supp f
\subset A$, $ \int_{\Omega}f \,d\mu = 0$ and $\mu \{t\dopu  f(t)=
-1\} > \mu(A) - \eps$. The collection of all  balanced
$\eps$-peaks on $A$ will be  denoted by $P(A,\eps)$.

\begin{thm}
$\narr(L_1) =  \{P(A,\eps)\dopu  A \in \Sigma,\ \eps > 0 \}^\sim$.
\end{thm}

\begin{proof}
Let $T \in \narr(L_1)$, $\delta, \eps > 0$, and $A \in \Sigma$.
Consider a slice in $ L_1$ of the form $$ S=\Bigl\{ f \in
B(L_1)\dopu   \int_{A}f \,d\mu > 1 - \delta \Bigr\}. $$ Applying
Lemma~\ref{Lnarrow} to this  slice, the elements
$x=-\chi_A/\mu(A)$, $ y=\chi_A/\mu(A)$ and $\delta$ we get a
function $v \in S$ such that
\begin{equation} \label{5}
\|v-\chi_A/\mu(A)\|>2-\delta, \quad
\|T(v-\chi_A/\mu(A))\|< \delta.
\end{equation}
Denote by $B$ the set $\{t\in A\dopu  v(t)>0\}$. The condition
$v \in S$ implies that $\|v-\chi_B v\| < \delta$, so
$$
\|v\chi_B-\chi_A/\mu(A)\|>2-2\delta.
$$
Next, introduce  $C=\{t\in A\dopu  v(t)>1/\mu(A)\}$. By the last inequality
$$
\|v\chi_C-\chi_A/\mu(A)\|>2-2\delta,\quad \|v-\chi_C v\| < 3\delta
$$
and
\begin{equation} \label{10}
\mu(C)< \delta \mu(A);
\end{equation}
to see this observe that 
\begin{align}
2-2\delta   &< 
\Bigl\| \chi_{B}v - \frac{\chi_{A}}{\mu(A)} \Bigr\| 
\le \int_{C} \Bigl( \chi_{B}v - \frac1{\mu(A)} \Bigr) \,d\mu
     + \frac1{\mu(A)} (\mu(A) - \mu(C) )  \nonumber\\
&\le 2 - 2 \frac{\mu(C)}{\mu(A)}. \nonumber
\end{align}
Put $f=(\mu(A)/\beta)\chi_C v - \chi_A$ with $\beta= \int_{C}v\, d\mu$ 
so that $\int_{\Omega} f \,d\mu =0$. 
Since $\int_{A} v \,d\mu> 1-\delta$ we have from  $\|v-\chi_{C}v\|<3\delta$
that $\beta\ge 1-4\delta$.
By (\ref{5}) we conclude that
$$
\|Tf\| = 
\mu(A) \Bigl\| T\Bigl( \frac{\chi_{C}v}\beta -
\frac{\chi_{A}}{\mu(A)} \Bigr) \Bigr\| \le
\mu(A) \Bigl( \|T\|\, \Bigl\| \frac{\chi_{C}v}\beta - v \Bigr\|+ \delta \Bigr) 
$$
and
$$
\Bigl\| \frac{\chi_{C}v}\beta - v \Bigr\| \le
\Bigl\| \frac{\chi_{C}v-v}\beta  \Bigr\| +
\Bigl\| \frac{v}\beta - v \Bigr\| \le
\frac{3\delta}\beta + \Bigl( \frac1\beta -1 \Bigr) \le
\frac{7\delta}{1-4\delta},
$$
and if $\delta$ is small enough, by (\ref{10})
$f \in P(A,\eps)$. This proves the inclusion
$\narr(L_{1}) \subset  \{P(A,\eps)\dopu  A \in \Sigma,\ \eps > 0 \}^\sim$.

To prove the opposite inclusion we use Proposition~\ref{Prop3.11}.
Let us fix
$T \in \{P(A,\eps)\dopu  A \in \Sigma,\ \eps > 0 \}^\sim$.
Let  $x,y\in S(L_{1})$, $y^*\in S(L_{\infty})$ and $\eps>0$ be
such that $\langle y^*,y\rangle > 1 - \eps$. Without loss of generality
we may assume that there is a partition $A_1, \dots , A_{n}$
of $\Omega$ such that the restrictions of $x$, $y$ and  $y^*$
on $A_k$ are constants, say $a_k$, $b_k$ and $c_k$ respectively.
By our assumption $T$ is unbounded from below on
each of  the $P(A_k,\delta)$ for every $\delta > 0$, $  k=1, \dots , n$.
Let us fix functions $f_k \in P(A_k,\delta)$ such that
$\|Tf_k\|  < \delta$, $  k=1, \dots , n$, and put
$$
v = \sum_{k=1}^n  b_k (\chi_{A_k} + f_k).
$$
By 
definition of balanced $\delta$-peaks $\langle y^*,v\rangle
> 1 - \eps$,
$\|v\| =  1$, and $\|T(y-v)\|$ and  $\mu(\supp v)$ become arbitrarily
small when $\delta$ is small enough. Thus $\delta$ can be chosen
so that $v$ fulfills the conditions $\|T(y-v)\|< \eps$ and
$\|x+v\| > 2 - \eps$.
\end{proof}

The characterisation of narrow operators on $L_1$
proved above looks similar to the definition of $L_1$-narrow
operators. It is easy to prove that every $L_1$-narrow
operator is narrow. We don't know whether the classes
of narrow operators and $L_1$-narrow operators on $L_1$ coincide.

The aim of the remainder
of this section is to construct an example of a strong
Daugavet operator on $L_{1}$
which is not narrow. In fact, we shall define a subspace
$Y\subset L_{1}[0,1]$ so that the quotient map $q\dopu L_{1}\to
L_{1}/Y$ is a strong Daugavet
operator, but $Y$ fails the Daugavet property. By
Theorem~\ref{theo5.3}, $q$ cannot be narrow. Likewise, $Y$ is almost
rich, but not rich.

Let $I_{n,k}=[\frac{k-1}{2^n}, \frac{k}{2^n})$ for $n\in \N_{0}$ and
$k=1,2,\dots,2^n$. Fix $N\in\N$. We define
\begin{eqnarray*}
g_{0,1} &=&
(2^N-1)\chi_{ I_{N,1} } - \chi_{ I_{0,1}\setminus I_{N,1} }   
\allowbreak     \\
g_{1,1} &=&
(2^{N^2-N}-1)\chi_{ I_{N^2,1} } - \chi_{ I_{N,1}\setminus I_{N^2,1} }   \\
g_{1,k} &=&  g_{1,1} (t- \textstyle \frac{k-1}{2^n} ),
\qquad k=2,\dots,2^N,\\
    &\vdots& \\
g_{n,1} &=&
(2^{N^{n+1}-N^n}-1)\chi_{ I_{N^{n+1},1} } - \chi_{ I_{N^n,1} \setminus
I_{N^{n+1},1}  }  \\
g_{n,k} &=&  g_{n,1} (t- \textstyle \frac{k-1}{2^{N^n}} ),
\qquad k=2,\dots,2^{N^n}.
\end{eqnarray*}
Denote by $P_{n}$ the ``peak set'' of the $n$'th generation, i.e.,
$$
P_{n} = \biggl\{t\in [0,1]\dopu \sum_{k=1}^{2^{N^n}} g_{n,k}(t)>0
\biggr\},
$$
and $P=\bigcup_{n} P_{n}$. Clearly $|P_{n}| = 2^{N^n}/ 2^{N^{n+1}} =
\bigl( 1/2^{N-1} \bigr)^{N^n}$  and $|P|\le 1/(2^N-1)$.
Notice also that $\int_{0}^1 g_{n,k}(t)\,dt =0$ for all $n$ and $k$.

First we formulate a lemma.
All the norms appearing below are $L_{1}$-norms.

\begin{lemma}\label{L4.4}
Let
$$
g= \sum_{n=0}^M \sum_{k=1}^{2^{N^n}} a_{n,k}g_{n,k}.
$$
Then
\[
\|g \chi_{[0,1] \setminus  P} \| \le 3 \|g \chi_{P}\|.
\]
\end{lemma}

\begin{proof}
 Denote
$$
g'' = \sum_{\supp g_{n,k}\subset P} a_{n,k}g_{n,k}, \quad g'=g-g''.
$$
Since $g'$ and $g$ coincide off $P$, we clearly have
\begin{equation}\label{eq6.2.1}
\|g'\chi_{[0,1] \setminus  P}\| = \|g\chi_{[0,1] \setminus  P}\|.
\end{equation}

We also have that
\begin{equation}\label{eq6.2.2}
\|g'\chi_{P}\| \le \|g\chi_{P}\|.
\end{equation}

Indeed,
we can write $P$ as a countable union of disjoint (half-open) intervals;
denote by $I$ any one of these. Then $g'$ is constant on $I$, and
$\int_{0}^1 g''(t)\,dt =0$. Hence
$$
\|g'\chi_{I}\| =
\biggl| \int_{0}^1 g'(t) \chi_{I}(t)\,dt \biggr| =
\biggl| \int_{0}^1 (g'(t) \chi_{I}(t) + g''(t) \chi_{I}(t))\,dt \biggr|
\le
\|g\chi_{I}\|.
$$
Summing up over all $I$ gives the result.

Next, we claim that
\begin{equation}\label{eq6.2.3}
\|g' \chi_{[0,1] \setminus  P} \| \le 3 \|g' \chi_{P}\|.
\end{equation}
To see this,
we label the intervals $I$ from the previous paragraph
as follows. For every $l\in\N$ write $B_{0}=P_{0}$ and $B_{l}=
P_{l}\setminus \bigcup_{i=1}^{l-1}P_{i}$. Each $B_{l}$ can be written
as $\bigcup_{d\in D_{l}} I_{N^{l+1},d}$ where $D_{l}$ is some subset
of $\{1,\dots, 2^{N^{l+1}} \}$  with cardinality ${ <2^{N^l}}$.
Let us write $g'= \sum_{n=0}^M \sum_{k=1}^{2^{N^n}} b_{n,k}g_{n,k}$.
We then have the estimates
$$
\int_{0}^1 |g'(t)\chi_{B_{0}}(t)|\,dt = |b_{0,1}| \frac{2^N-1}{2^N}
$$
and
\begin{eqnarray*}
\int_{0}^1 |g'(t)\chi_{B_{l}}(t)|\,dt
&=&
\sum_{d\in D_{l}} \int_{I_{N^{l+1},d}} \biggl| -b_{0,1}
   - \sum_{n=1}^{l-1} \sum_{k=1}^{2^{N^n}} b_{n,k} \chi_{\supp g_{n,k}}
 \\
&& \qquad \mbox{}
      + b_{l, (d-1)/(2^{N-1})^{N^l} +1} \bigl(  2^{N^{l+1} - N^l} -1
            \bigr) \biggr| \,dt  \\
&\ge&
\sum_{k=1}^{2^{N^l}} \Bigl( \frac1{2^{N^l}} - \frac1{2^{N^{l+1}}}
     \Bigr) |b_{l,k}|
\\
&& \qquad \mbox{}
- \frac1{(2^{N-1})^{N^l}}       |b_{0,1}| -
\frac1{(2^{N-1})^{N^l}} \sum_{n=1} ^{l-1} \sum_{k=1}^{2^{N^l}}
|b_{n,k}|.
\end{eqnarray*}
Summing up over all $l$ gives us
\begin{eqnarray*}
\int_{P} |g'(t)|\,dt  &\ge&
|b_{0,1}| \biggl( \frac{2^N-1}{2^N} - \sum_{m=1}^\infty
\frac1{(2^{N-1})^{N^m}} \biggr)
\\
&& \qquad \mbox{}
+ \sum_{l=1}^\infty \biggl( \frac1{2^{N^l}} - \frac1{2^{N^{l+1}}} -
   \sum_{m=l+1}^\infty \frac1{(2^{N-1})^{N^m}}  \biggr)
    \sum_{k=1}^{2^{N^l}} |b_{l,k}| \\
&\ge&
\frac12 |b_{0,1}|
+ \frac12 \sum_{l=1}^\infty \frac1{2^{N^l}}
    \sum_{k=1}^{2^{N^l}} |b_{l,k}| .
\end{eqnarray*}
On the other hand, by the triangle inequality
$$
\int_{0}^1 |g'(t)|\,dt \le 2 \biggl( |b_{0,1}| +
    \sum_{l=1}^\infty \frac1{2^{N^l}}
    \sum_{k=1}^{2^{N^l}} |b_{l,k}| \biggr),
$$
hence the claim follows.

The lemma now results from (\ref{eq6.2.1})--(\ref{eq6.2.3})
\end{proof}

\begin{theo} \label{theo6.3}
Let $Y_{N}\subset L_{1}[0,1]$ be the closed subspace 
generated by the system $ \{g_{n,k}\}$ and the constants.
Then the quotient map
$q_{N}\dopu L_{1}\to L_{1}/Y_{N}$ is a strong Daugavet operator
for all $N$, but
$Y_{N}$ fails the Daugavet property if $N\ge4$.
\end{theo}

\begin{proof}
Let us fix $x,y\in S(L_{1})$ and $\eps>0$. Without loss of generality
we may assume that $x= \sum_{k=1}^{2^{N^n}} a_{n,k} \chi_{I_{n,k}}$ for
a big enough $n$ to be chosen later.

Put $h= \sum_{k=1}^{2^{N^n}} a_{n,k} g_{n,k}$. Then
$$
x+h =  \sum_{k=1}^{2^{N^n}} 2^{ N^{n+1}-N^n } \chi_{N^ {n+1},d_{n,k}}
a_{n,k}
$$
with $d_{n,k}= 1+(k-1) (2^{N-1})^{N^n}$. So
$$
\|x+h\| = \sum_{k=1}^{2^{N^n}} \frac{|a_{n,k}|}{2^{N^n}} = \|x\| =1,
$$
and $\supp (x+h)\subset P_{n}$. Since $|P_{n}|\to0$ we can pick $n$
big enough to satisfy $\|x+h+y\|>2-\eps$.
This shows that  $q_{N}$ is a strong Daugavet operator.

To show that
$Y_{N}$ fails the Daugavet property if $N\ge4$,
take $g^*=\chi_{[0,1] \setminus  P}\in Y_{N}^*$  and $\eps=2|P|$. Since
$\eins\in S(Y_{N})$, we get
$$
\|g^*\| \ge g^*(\eins) = 1-\eps/2 >1-\eps.
$$
Thus, $S({g^*,\eps})\cap B(Y_{N})\neq \emptyset$. We show that there
is no $f$ in this slice such that $\|f-\eins\|>2-\eps$.

Suppose, on the contrary, that there is such an $f$. Without loss of
generality we can assume that
$$
f=a_{0}\eins + g
$$
where $g$ is as in Lemma~\ref{L4.4}.

It follows from our conditions that
\begin{equation}\label{eq1}
\|f\chi_{P}\| = \int_{P} |f(t)|\,dt = \|f\| - g^*(|f|) \le
1-g^*(f) < \eps.
\end{equation}
Hence,
$$
1\ge \int_{0}^1 f(t)\,dt = \int_{P} f(t)\,dt + g^*(f) > 1-2\eps,
$$
and since $\int_{0}^1 f(t)\,dt =a_{0}$, we get
\begin{equation}\label{eq2}
1-2\eps<a_{0}\le1.
\end{equation}
By (\ref{eq1}) and (\ref{eq2}),
\begin{equation}\label{eq3}
\|g\chi_{P}\| \le \eps + |P| < 2\eps,
\end{equation}
thus (\ref{eq2}) and (\ref{eq3}) yield
$$
\|g\chi_{[0,1] \setminus  P}\|
\ge \|g\|-2\eps
= \|f- a_{0}\eins\| -2\eps
 \ge \|f-\eins\|-4\eps > 2-5\eps.
$$
But now Lemma~\ref{L4.4} and (\ref{eq3}) imply
$$
2-5\eps<   \|g\chi_{[0,1] \setminus  P}\| \le 3\|g\chi_{P}\| <6\eps,
$$
which yields $\eps>2/11$, i.e., $|P|>1/11$, which is false for $N\ge4$.
\end{proof}

Theorems~\ref{theo6.3} and~\ref{theo5.3} immediately yield
the following result.

\begin{thm} \label{theo6.7}
There is an almost rich subspace  of $L_{1}[0,1]$
which fails the Daugavet property and hence fails to be rich.
Thus, on $L_{1}[0,1]$ the class of strong
Daugavet operators does not coincide with the class of narrow
operators.
\end{thm}


\section{Questions}

We finish this paper with some questions which have remained open. We
intend to deal with these problems in a future publication.

\begin{enumerate}
\item \label{qu1}
Does the class of narrow operators on a Banach space $X$ form a subsemigroup
of $\Op(X)$? [Added in the final version: We have recently
constructed a counterexample on the space $X=C([0,1], L_{1}[0,1])$.]
\item \label{qu2}
Is every narrow operator on $L_{1}$ also $L_{1}$-narrow?
\item \label{qu3}
Is the sum of two $L_{1}$-narrow operators from $L_{1}$  to $L_{1}$ again
$L_{1}$-narrow? This question is clearly related to the previous
ones; we remark that the proof in \cite[p.~69]{PliPop} which
purportedly shows this to be true appears to have a gap.
\item \label{qu4}
If $X$ has the Daugavet property, does $X$ have a subspace isomorphic
to $\ell_{2}$?
\item \label{qu5}
If $T$ is an operator on a space $X$ with the Daugavet property which
does not fix a copy of $\ell_{2}$, is $T$ then narrow? We remark that
the answer is affirmative in the case $X=C[0,1]$ by our
Theorem~\ref{l1} and a result due to Bourgain \cite{Bour-JOT}.
\end{enumerate}



\end{document}